\documentclass[12pt,reqno]{amsart}

\usepackage[utf8]{inputenc} %input encoding for non-ASCII characters
\usepackage[T1]{fontenc} %T1 font encoding

\usepackage{amscd,graphicx}
\usepackage[all]{xy}
\usepackage{fullpage}
\usepackage{subsupscripts}

\usepackage{ifthen}

\newboolean{workmode}
\setboolean{workmode}{false}

\newboolean{sections}
\setboolean{sections}{true}

\usepackage[indent=0.5cm]{parskip}

\usepackage{multicol}
	\AtBeginDocument{\addtocontents{toc}{\protect\setlength{\parskip}{0pt}}}
	\makeatletter
		\renewcommand*\l@subsection{\@tocline{2}{0pt}{2.75pc}{5pc}{}}
	\makeatother

	\usepackage[backend=bibtex,style=alphabetic,defernumbers=true,giveninits=false,doi=false,isbn=false,url=false,useprefix=true,maxnames=10,
	minnames=9,maxbibnames=99,maxalphanames=10,minalphanames=9,sorting = nyt,backref=false]{biblatex} %bibliography style
	\bibliography{bibliography.bib} %selects the correct bibliography file
	\DeclareFieldFormat{labelalpha}{\mkbibbold{#1}}
	\DeclareFieldFormat{extraalpha}{\mkbibbold{\mknumalph{#1}}}
	\AtBeginBibliography{
		\DeclareFieldFormat{labelalpha}{#1}
		\DeclareFieldFormat{extraalpha}{\mknumalph{#1}}
	}

\usepackage{amsmath}

	\DeclareFontFamily{U}{mathx}{}
	\DeclareFontShape{U}{mathx}{m}{n}{<-> mathx10}{}
	\DeclareSymbolFont{mathx}{U}{mathx}{m}{n}
	\DeclareMathAccent{\widehat}{0}{mathx}{"70}
	\DeclareMathAccent{\widecheck}{0}{mathx}{"71}

	\DeclareFontFamily{U}{mathx}{\hyphenchar\font45}
    \DeclareFontShape{U}{mathx}{m}{n}{<-> mathx10}{}
    \DeclareSymbolFont{mathx}{U}{mathx}{m}{n}
    \DeclareMathAccent{\widebar}{0}{mathx}{"73}    
	
\usepackage{amssymb}

\usepackage{amsthm}
\usepackage{thmtools} %fixes cref names
\usepackage{tikz-cd} %commutative diagrams
\usepackage{rotating} %for rotating symbols
\usepackage{xcolor}
\definecolor{jaw}{rgb}{0,.5,0}  % for Jordan's comments
\definecolor{forestgreen}{rgb}{.2,.6,.2} % for citations
\definecolor{darkturquoise}{RGB}{33,108,115} %for links

\usepackage[colorlinks=true,linkcolor={darkturquoise},citecolor={forestgreen},urlcolor={darkturquoise}]{hyperref}
\usepackage[capitalise,nameinlink,noabbrev]{cleveref}
	\crefname{subsection}{Subsection}{subsections}
	\crefdefaultlabelformat{#2\textbf{#1}#3} %makes the numbers bold

\usepackage{fontawesome}
\newcommand{\ifwork}[1]{\ifthenelse{\boolean{workmode}}{#1}{}}
\newcommand{\comment}[1]{}
\newcommand{\mute}[1]{}
\newcommand{\printname}[1]{}

\ifwork{
	\renewcommand{\comment}[1]{{\marginpar{*}\ \scriptsize{#1}\ }}
	\renewcommand{\printname}[1]
	{{\color{brown}{\makebox[0pt]{\hspace{-1in}\raisebox{8pt}{\tiny #1}}}}}
}

\newcommand{\labell}[1]{\label{#1} \printname{#1}}

		\newcommand{\calA}{{\mathcal{A}}}
		\newcommand{\calB}{{\mathcal{B}}}
		\newcommand{\calC}{{\mathcal{C}}}
		\newcommand{\calD}{{\mathcal{D}}}
		
		\newcommand{\calF}{{\mathcal{F}}}
		\newcommand{\calG}{{\mathcal{G}}}
		\newcommand{\calH}{{\mathcal{H}}}
		\newcommand{\calI}{{\mathcal{I}}}
		\newcommand{\calJ}{{\mathcal{J}}}
		\newcommand{\calK}{{\mathcal{K}}}
		\newcommand{\calL}{{\mathcal{L}}}
		\newcommand{\calM}{{\mathcal{M}}}
		\newcommand{\calN}{{\mathcal{N}}}

		\newcommand{\calU}{{\mathcal{U}}}
		\newcommand{\calV}{{\mathcal{V}}}
		\newcommand{\calW}{{\mathcal{W}}}
		\newcommand{\calX}{{\mathcal{X}}}

		\newcommand{\NN}{{\mathbb{N}}}  %natural numbers
		\newcommand{\RR}{{\mathbb{R}}}  %real numbers or cartesian space
		\renewcommand{\SS}{{\mathbb{S}}} %sphere
		\newcommand{\ZZ}{{\mathbb{Z}}}  %integers

		\newcommand{\ana}{{\operatorname{ana}}} %anafunctor
		\newcommand{\bi}{{\operatorname{bi}}} %bi
		\newcommand{\Bicat}{{\mathbf{Bicat}}} %weak categorical thing of bicategories
		\newcommand{\DBiBund}{{\mathbf{DBiBund}}} % diffeological bibundles
		\newcommand{\DGpoid}{{\mathbf{DGpoid}}} % diffeological groupoids
		\newcommand{\Diffeol}{\mathbf{Diffeol}} %diffeological spaces
		\newcommand{\ev}{{\operatorname{ev}}}   %evaluation map
		\newcommand{\fibr}{{\operatorname{fibr}}} % fibration map
		\newcommand{\Hom}{{\operatorname{Hom}}} %Hom
		\newcommand{\id}{{\operatorname{id}}}  %identity
		\newcommand{\inv}{{\operatorname{inv}}} %inversion
		\newcommand{\Lie}{{\operatorname{Lie}}} %Lie
		\newcommand{\LieGpoid}{{\mathbf{LieGpoid}}} %Lie groupoid
		\newcommand{\mult}{{\operatorname{m}}} % groupoid multiplication
		\newcommand{\Neb}{{\operatorname{Neb}}} %nebula
		\newcommand{\pr}{{\operatorname{pr}}} %projection
		\newcommand{\PR}{{\operatorname{Pr}}} %Projection
		\newcommand{\src}{{\operatorname{s}}} % source
		\newcommand{\trg}{{\operatorname{t}}} % target
		\newcommand{\unit}{{\operatorname{u}}} % unit
		\newcommand{\arr}[1]{{\,\raisebox{1mm}{$\overset{\scriptstyle{#1}}{\curvearrowright}$}\,}} %groupoid arrow
		\newcommand{\del}{\partial}  % del
		\newcommand{\eps}{\varepsilon}  %var epsilon
		\newcommand{\ftimes}[2]{{\lrsubscripts{\times}{#1}{#2}}} %fibre product
		\newcommand{\toto}{{~\rightrightarrows~}} %groupoid double arrow
		\newcommand{\LAalong}[1]{{\,\curvearrowright\hspace{-5pt}^{#1}\hspace{1pt}}} % left groupoid action
		\newcommand{\onto}{{\,\scalebox{2}[1]{$\twoheadrightarrow$}\,}} %surjective arrow
		\newcommand{\ot}[1]{{\underset{\raisebox{1mm}{$\quad\scriptstyle{#1}$}}{\overset{\raisebox{-1mm}{$\quad\scriptstyle{\simeq}$}}{\,\scalebox{2}[1]{$\twoheadleftarrow$}}}\,}} %left arrow
		\newcommand{\RAalong}[1]{{\hspace{-1pt}~^{#1}\hspace{-5pt}\curvearrowleft\,}} % right groupoid action
		\newcommand{\subdwe}{{\,\overset{\simeq}{\scalebox{2}[1]{$\twoheadrightarrow$}}\,}} %subductive weak equivalence
		\newcommand{\uto}[1]{{\,\underset{\raisebox{1mm}{$\scriptstyle{#1}$}}{\longrightarrow}\,}} %right arrow
		\newcommand{\we}{{\,\overset{\simeq}{\longrightarrow}\,}} %weak equivalence
		\newcommand{\wtimes}[2]{{\lrsubscripts{\overset{\raisebox{-2pt}{\tiny$\mathrm{w}$}}{\times}}{#1}{#2}}} %weak pullback
		\usepackage{faktor}

		\makeatletter
		\DeclareRobustCommand*{\mfaktor}[3][]
		{
		   { \mathpalette{\mfaktor@impl@}{{#1}{#2}{#3}} }
		}
		\newcommand*{\mfaktor@impl@}[2]{\mfaktor@impl#1#2}
		\newcommand*{\mfaktor@impl}[4]{
		   \settoheight{\faktor@zaehlerhoehe}{\ensuremath{#1#2{#3}}}%
		   \settoheight{\faktor@nennerhoehe}{\ensuremath{#1#2{#4}}}%
		      \raisebox{-0.5\faktor@zaehlerhoehe}{\ensuremath{#1#2{#3}}}%
		      \mkern-4mu\diagdown\mkern-5mu%
		      \raisebox{0.5\faktor@nennerhoehe}{\ensuremath{#1#2{#4}}}%
		}
		\makeatother

	\newcommand{\ifsection}[2]{\ifthenelse{\boolean{sections}}{#1}{#2}} % Sets numbering for environments based on sections boolean value (at top of this document).

		\theoremstyle{plain}
		\ifsection{
			\newtheorem{theorem}{Theorem}[section] % follows section numbers
		}
		{
			\newtheorem{theorem}{Theorem} % does not follow section numbers
		}
		\newtheorem*{theorem*}{Theorem}
		\newtheorem{proposition}[theorem]{Proposition}
		\newtheorem{corollary}[theorem]{Corollary}
		\newtheorem*{conjecture*}{Conjecture}
		
		\newtheorem{lemma}[theorem]{Lemma}
		\newtheorem*{lemma*}{Lemma}

		\theoremstyle{definition}
		\newtheorem{definitionx}[theorem]{Definition}
		\newtheorem{examplex}[theorem]{Example}
		\newtheorem{examplesx}[theorem]{Examples}

		\newtheorem{remarkx}[theorem]{Remark}

	\newenvironment{definition}
		{\pushQED{\qed}\definitionx}
		{\popQED\enddefinitionx}

	\newenvironment{example}
		{\pushQED{\qed}\examplex}
		{\popQED\endexamplex}

	\newenvironment{examples}
		{\pushQED{\qed}\examplesx}
		{\popQED\endexamplesx}

	\newenvironment{remark}
		{\pushQED{\qed}\remarkx}
		{\popQED\endremarkx} % All non-paper-specific macros, theorem styles, work mode and section options, etc. are in the file ``metamatter.tex''.

	\subjclass[2020]{Primary 22A22; Secondary 58H05, 57R55} 

	\keywords{diffeology, groupoid, anafunctor, Morita equivalence}

\begin{document}

	\author{Jordan Watts}
	\address{Department of Mathematics, Central Michigan University, Pearce Hall 214, Mount Pleasant, MI 48859, USA}
	\email{jordan.watts@cmich.edu}
	\title{The Bicategory of Lie Groupoids within Diffeological Groupoids}
	\date{\today}

	\begin{abstract}
		We consider the localisation of the $2$-category of diffeological groupoids at weak equivalences from the perspective of anafunctors, and with this language, prove that the localisation of the $2$-category of Lie groupoids is an essentially full sub-bicategory of that of diffeological groupoids.  In particular, we solve the open problem affirmatively of whether two Lie groupoids that are diffeologically Morita equivalent are Morita equivalent in the usual Lie sense.  
	\end{abstract}

	\maketitle

	\thispagestyle{empty}

%%%%%%%%%%%%%%%%%%%%%%%%%%%%%%%%%%%%%%%%%%%%%%%%
\section{Introduction}\labell{s:intro}
%%%%%%%%%%%%%%%%%%%%%%%%%%%%%%%%%%%%%%%%%%%%%%%%

	Diffeological groupoids and the application of diffeology to the study of Lie groupoids and Lie algebroids are current and important trends in geometry.  They appear in work on general relativity \cite{BFW}, in which a diffeological groupoid describes the choices of embeddings of an initial space-like hypersurface in a lorentzian spacetime, up to a prescribed equivalence.  They show up in the study of singular subalgebroids of Lie algebroids \cite{AZ}, the integration of Lie algebroids \cite{villatoro:integrability}, and the holonomy and fundamental groupoids of a singular foliation \cite{GV}.  Diffeological groupoids have become crucial in the study of the ``higher geometric'' version of loop spaces: loop stacks; in \cite{RV:loop}, the authors show that the stack $\Hom(\SS^1,\calX)$ is presentable by not just a diffeological groupoid, but a Fr\'echet-Lie groupoid, where $\calX$ is a differentiable stack. As Fr\'echet manifolds have been shown to form a full subcategory of the category of diffeological spaces \cite[Section 3]{losik}, here is an example (of many) where working in the diffeological category is very beneficial to infinite-dimensional differential geometry. In \cite{RV:orbifold}, the authors announce that they have extended the results of \cite{RV:loop} to stacks $\Hom(M,\calX)$ where $M$ is a compact manifold.  Yet another place diffeological groupoids appear is as inertia groupoids, which play an important role in the K-theory and Chen-Ruan cohomology for orbifolds \cite{ALR}.

	Thus it became natural for a rigorous foundation for diffeological groupoids and their Morita equivalence to be developed.  The thesis and subsequent paper of van der Schaaf \cite{vdS:thesis,vdS:morita} provide such a foundation in terms of bibundles, extending the theory of Lie groupoids and their bibundles to the diffeological realm; indeed, the bicategory of Lie groupoids, right principal bibundles, and bi-equivariant diffeomorphisms forms a sub-bicategory of the diffeological version.  An important open question \cite[Question 7.6]{vdS:morita} is whether a Morita equivalence between two Lie groupoids \emph{in the diffeological bicategory} is, in fact, a Lie Morita equivalence; that is, whether the biprincipal bibundle between the two Lie groupoids representing the Morita equivalence is actually a bibundle from the Lie bicategory (\emph{i.e.}\ a smooth manifold).  An affirmative answer is one of the motivations for the present manuscript.
	
	In this paper, we construct a localisation of diffeological groupoids at weak equivalences using the anafunctor (or $J$-fraction) setting of Roberts \cite{roberts:anafunctors}.  Other options in which to construct a localisation include the generalised morphisms of Pronk and Pronk-Scull \cite{pronk,pronk-scull:bicateg}, the bibundle formalism of van der Schaaf mentioned above, or stacks.  The work of Pronk and Pronk-Scull works in great generality, and as such, computations (especially those involving $2$-cells) can get very complicated.  On the other hand, bibundles are very rigid, as their definition pins down the exact geometric attributes required to invert weak equivalences.  The language of stacks, as categories fibred in groupoids, can very quickly lose any sense of the geometry at play.  Thus it is preferred by this author to utilise a happy medium that seems to work best for the purposes at hand.  Anafunctors originally were introduced by Makkai \cite{makkai} in order to allow one to discuss $1$-cells in a $2$-category without using the axiom of choice.  Bartels \cite{bartels} develops the theory further in terms of internal categories.  Roberts emphasises the use of coverings to achieve the localisation \cite{roberts:internal,roberts:anafunctors}, which adapts well to the diffeological setting.
	
	As alluded to above, in order to localise diffeological groupoids at weak equivalences in terms of anafunctors, certain prerequisites need to be met; this is the topic of \cref{s:2-site}.  A description of the resulting bicategory of diffeological groupoids and anafunctors is found in \cref{s:anafunctors}.  We move to Lie groupoids in \cref{s:lie gpds}, where we prove that the localised bicategory of Lie groupoids forms an essentially full sub-bicategory of the localised bicategory of diffeological groupoids.  While we do not require the arrow spaces of the Lie groupoids to be Hausdorff nor second-countable, restricting our attention to those that do satisfy these conditions in turn lead to essentially full sub-bicategories. \cref{s:bibundles} reviews the diffeological bibundle setting of van der Schaaf, and proves that this bicategory is also a localisation of the $2$-category of diffeological groupoids at weak equivalences, and hence equivalent to the bicategory constructed using anafunctors.  The open question of van der Schaaf is hence answered.  We end the paper with some examples and applications in \cref{s:apps and examples}.
	
	The preliminaries for this paper are potentially substantial, and we provide a preliminary section \cref{s:preliminaries} more to set notation than to review material, although we provide plenty of examples.  Therefore, we recommend the reader refer to the following sources as needed: \cite{iglesias} for details on diffeology, \cite{JY} for details on $2$-categories and bicategories, and \cite{vdS:morita} for details on the bicategory of diffeological groupoids with bibundles.  We try to pinpoint the relative places within these as they are used.  
	
	\noindent\textbf{Acknowledgements:} The author wishes to thank Dorette Pronk and Laura Scull for many conversations on subtleties of bicategorical theory and localisation, and Nesta van der Schaaf for his help as this project was just getting off the ground.  We are also indebted to Laura Scull for suggestions on notation, which helped make the paper much easier to read, and David Michael Roberts for his feedback on a preprint out of which this paper arose.  Finally, we thank the anonymous referee for excellent suggestions, improving the quality of this paper.

%%%%%%%%%%%%%%%%%%%%%%%%%%%%%%%%%%%%%%%%%%%%%%%%
\section{Preliminaries}\labell{s:preliminaries}
%%%%%%%%%%%%%%%%%%%%%%%%%%%%%%%%%%%%%%%%%%%%%%%%

	We begin by setting our notation and giving a quick overview of diffeology and groupoids.  For more details on diffeology, see \cite{iglesias}. Throughout this paper, all manifolds are smooth, second-countable, and Hausdorff, unless stated otherwise.

	\begin{definition}[Diffeology]\labell{d:diffeology}
        	Let $X$ be a set.  A \textbf{cartesian open set} is an open set of some cartesian space $\RR^n$ for $n\in\NN=\{0,1,\dots\}$.  A \textbf{parametrisation} of $X$ is a set-theoretical function $p\colon U\to X$ where $U$ is a cartesian open set.  Given a parametrisation $p$, we denote by $U_p$ the domain of $p$, and for convenience we will often utilise the notation $p\colon u\mapsto x_u$.  (Without loss of generality we will typically be able to assume that $0\in U_p$, and so utilise $x_0$ in this context without justification.  See, for instance, the notation for a local subduction in \cref{d:ind subd}.) A \textbf{diffeology} $\calD$ on $X$ is a family of parametrisations satisfying:
	        \begin{enumerate}
	            \item \textbf{(Covering Axiom)} All constant parametrisations are contained in $\calD$;
	            \item \textbf{(Locality Axiom)} If $p$ is a parametrisation of $X$ and $U_p$ admits an open cover $\{U_\alpha\}$ for which $p|_{U_\alpha}\in\calD$ for each $\alpha$, then $p\in\calD$;
	            \item \textbf{(Smooth Compatibility Axiom)} If $p\in\calD$ and $f\colon V\to U_p$ is a smooth map from a cartesian open set $V$, then $p\circ f\in\calD$.
	        \end{enumerate}
	        We refer to $(X,\calD)$ as a \textbf{diffeological space}, and the parametrisations in $\calD$ as \textbf{plots}.  We will typically drop the notation $\calD$, and denote the diffeology of a diffeological space $X$ as $\calD_X$ when needed.

        	Given two diffeological spaces $X$ and $Y$, a set-theoretical function $F\colon X\to Y$ is \textbf{(diffeologically) smooth} if $F\circ p\in\calD_Y$ for each plot $p\in\calD_X$.  A smooth bijection whose inverse is also smooth is a \textbf{diffeomorphism}.
	\end{definition}

	The category $\Diffeol$ of diffeological spaces is a complete cocomplete quasi-topos \cite{BH}; in particular, it admits subobjects, quotients, products, coproducts, and function spaces.  
    
	\begin{examples}[Examples of Diffeological Spaces]\labell{x:diffeology}
    		\noindent
    		\begin{enumerate}
	    		\item\labell{i:manifolds}Given a (smooth) manifold $M$, the \textbf{standard manifold diffeology} on $M$ is the collection of all smooth parametrisations of $M$ in the classical sense.  In fact, the category of smooth manifolds can be identified with a full subcategory of $\Diffeol$.
	    		
	    		\item\labell{i:gen family}Fix a diffeological space $X$.  A \textbf{(covering) generating family} of the diffeology $\calD_X$ is a family of plots $\calF\subseteq\calD_X$ satisfying: for each $p\in\calD_X$ there exist an open cover $\{U_\alpha\}$ of $U_p$ and for each $\alpha$, a plot $q_\alpha\colon V_\alpha\to X$ in $\calF$ and a smooth function $f_\alpha\colon U_\alpha\to V_\alpha$ such that $p|_{U_\alpha}=q_\alpha\circ f_\alpha$.  The \textbf{nebula} $\Neb(\calF)$ of $\calF$ is the coproduct $\coprod_{q\in\calF} U_q$, which comes equipped with the (smooth) \textbf{evaluation map} $\ev_\calF\colon\Neb(\calF)\to X$ sending $u\in U_q$ to $q(u)$ for each $q\in\calF$.
	    		
	    		\item\labell{i:subset}Given a diffeological space $X$ and a subset $Y\subseteq X$, the \textbf{subset diffeology} on $Y$ is the subset of all plots of $X$ with image in $Y$.
	    		\item\labell{i:product}Given diffeological spaces $X$ and $Y$, the \textbf{product diffeology} on $X\times Y$ is the collection of parametrisations $(p,q)\colon U\to X\times Y$ for which $p\in\calD_X$ and $q\in\calD_Y$.
	    		\item\labell{i:quotient}Given a diffeological space $X$ and an equivalence relation $\sim$ on $X$ with quotient map $\pi$, the \textbf{quotient diffeology} on $X/\!\sim$ is the collection of all parametrisations $p$ for which $U_p$ admits an open cover $\{U_\alpha\}$ and for each $\alpha$ a plot $q_\alpha\colon U_\alpha\to X$ such that $p|_{U_\alpha}=\pi\circ q_\alpha$.\qedhere
		\end{enumerate}
	\end{examples}
	
	We will need to consider special types of smooth maps between diffeological spaces.

	\begin{definition}[Inductions \& Subductions]\labell{d:ind subd}
        	A smooth injection $F\colon X\to Y$ between diffeological spaces is an \textbf{induction} if every plot $q$ of $Y$ contained in the image of $F$ is equal to $F\circ p$ for some plot $p$ of $X$.
    
	        A smooth surjection $F\colon X\to Y$ between diffeological spaces is a \textbf{subduction} if for every plot $p$ of $Y$, there is an open cover $\{U_\alpha\}$ of $U_p$, and for each $\alpha$, a plot $q_\alpha\colon U_\alpha\to X$, called a \textbf{(local) lift against $F$}, such that $$p|_{U_\alpha}=F\circ q_\alpha.$$  If further, for every $u\in U_p$ and $x\in F^{-1}(p(u))$, there is an open neighbourhood $V$ of $u$ and a lift $q\colon V\to X$ with $q(u)=x$, then $F$ is a \textbf{local subduction}; in this case we call $q$ a \textbf{(local) lift against $F$ through $x$}.  Notationally, if $p\colon u\mapsto y_u$ is a plot, we typically can without loss of generality assume that $0\in U_p$; now if $x_0\in F^{-1}(y_0)$, then we will denote the lift $q$ by $q\colon u\mapsto x_u$, which implies that it is a lift through $x_0$.
	\end{definition}

	In practice when working with subductions, we often just shrink the domain of a plot $p$ and search for a global lift $q$.

	\begin{examples}[Examples of Inductions \& Subductions]\labell{x:ind subd}
	\noindent
        	\begin{enumerate}
        		\item\labell{i:inclusion}Given a diffeological space $X$ and a subset $Y\subseteq X$, the inclusion map $Y\hookrightarrow X$ is an induction provided we equip $Y$ with the subset diffeology.
        	
            		\item\labell{i:embd submfld}An embedded submanifold of a manifold $i\colon M\to N$ is an induction.
  
            		\item\labell{i:nebula}Given a diffeological space $X$ and a generating family $\calF$ of its diffeology, the evaluation map $\ev_\calF\colon\Neb(\calF)\to X$ is a subduction.
            		\item\labell{i:quotient subd}Given a diffeological space $X$ with an equivalence relation $\sim$, the quotient map $X\to X/\!\sim$ is a subduction.

            		\item\labell{i:not local subd}Let $M$ be a manifold of dimension greater than $0$, let $x_0\in M$ be fixed, and let $F\colon M\amalg M\to M$ be given by $F(x)=x_0$ for all $x$ in the first copy of $M$, and $F(x)=x$ for all $x$ in the second copy of $M$.  Then $F$ is a subduction, but not a local subduction.

            		\item\labell{i:submersion}A smooth map $F\colon M\to N$ between manifolds is a local subduction if and only if it is a surjective submersion. 
            		
            		\item\labell{i:diffeo inj subd}A smooth map $F\colon X\to Y$ between diffeological spaces is a diffeomorphism if and only if it is an injective subduction.
            
            		\item\labell{i:comp ind subd}The composition of two inductions is an induction, and the composition of two subductions is a subduction.
            		
            		\item\labell{i:pullback subd}Given the following pullback diagram of diffeological spaces, if $f$ is a subduction, then so is $\pr_2$ (and symmetrically, if $g$ is a subduction, then so is $\pr_1$):
            			$$\begin{gathered}[b]\xymatrix{
            				X\ftimes{f}{g}Y \ar[r]^{\quad\pr_2} \ar[d]_{\pr_1} & Y \ar[d]^{g} \\
            				X \ar[r]_{f} & Z. \\
				}\\[-\dp\strutbox]\end{gathered}\eqno\qedhere$$
       		\end{enumerate}
    	\end{examples}

	We are ready for the definition of a diffeological groupoid.

	\begin{definition}[Diffeological Groupoid]\labell{d:diffeol gpd}
        	A \textbf{diffeological groupoid} $\calG = \calG_1\toto\calG_0$ is a small groupoid whose sets of objects $\calG_0$ and arrows $\calG_1$ are diffeological spaces for which the following structure maps are smooth:
	        \begin{enumerate}
	            	\item The \textbf{source map} $\displaystyle\src_\calG\colon\calG_1\to\calG_0\colon x\arr{g}y~\mapsto~ x,$
	            	\item The \textbf{target map} $\displaystyle\trg_\calG\colon\calG_1\to\calG_0\colon x\arr{g}y~\mapsto~ y,$
	            	\item The \textbf{unit map} $\displaystyle\unit_\calG\colon\calG_0\to\calG_1\colon x\mapsto u_x,$
	            	\item The \textbf{multiplication map} $\displaystyle\mult_\calG\colon\calG_1\ftimes{\src}{\trg}\calG_1\to\calG_1\colon (g,h)\mapsto gh,$
	            	\item The \textbf{inversion map} $\displaystyle\inv_\calG\colon\calG_1\to\calG_1\colon x\arr{g}y~\mapsto~ y\arr{g^{-1}}x.$
	        \end{enumerate}
We will drop the subscripts from the structure maps above when the notation becomes too cluttered.

        	Given diffeological groupoids $\calG$ and $\calH$, a functor $\varphi\colon\calG\to\calH$ is \textbf{smooth} if the map between arrows $\varphi_1\colon\calG_1\to\calH_1$ is smooth.  A smooth functor admitting a smooth inverse functor is an \textbf{isomorphism} of diffeological groupoids.

        	Given smooth functors $\varphi,\varphi'\colon\calG\to\calH$, a natural transformation $S\colon \varphi\Rightarrow \varphi'$ is \textbf{smooth} if the underlying map $S\colon\calG_0\to\calH_1\colon x\mapsto S_x$ is smooth.

        	Diffeological groupoids with smooth functors and smooth natural transformations form a strict $2$-category, denoted $\DGpoid$.
	\end{definition}
	
	\begin{remark}\labell{r:diffeol gpd}
		It follows from the definition of a diffeological groupoid that the source and target maps are automatically subductions,  the unit map an induction, and inversion a diffeomorphism.  Given a smooth functor $\varphi\colon\calG\to\calH$, the map on objects $\varphi_0\colon\calG_0\to\calH_0$ is equal to $s_\calH\circ\varphi_1\circ\unit_\calG$, and hence is automatically smooth and determined completely by $\varphi_1$.
	\end{remark}

	\begin{examples}[Examples of Diffeological Groupoids]\labell{x:diffeol gpd}
		\noindent
		\begin{enumerate}
			\item\labell{i:lie gpd}A Lie groupoid is an example of a diffeological groupoid, using the standard diffeological structures on the spaces of objects and arrows.  By definition, the source and target maps are required to be submersions (\emph{i.e.} local subductions), which enable the multiplication map to have a manifold for its domain.
			
			\item\labell{i:not Lie}\cref{i:not local subd} of \cref{x:ind subd} provides an example of a diffeological groupoid that is not Lie, even though the object and arrow spaces are manifolds and all structure maps smooth ($F$ is both source and target).
						
			\item\labell{i:product gpd}Given diffeological groupoids $\calG$ and $\calH$, the \textbf{product groupoid} is the diffeological groupoid $\calG\times\calH$, where the object and arrow spaces are the products of the corresponding spaces of $\calG$ and $\calH$, and all structure maps are the natural ones induced by the product.
			
			\item\labell{i:pullback gpd}Let $\calG$, $\calH$, and $\calK$ be diffeological groupoids, and let $\varphi\colon\calG\to\calK$ and $\psi\colon\calH\to\calK$ be smooth functors.  The \textbf{(strict) pullback groupoid} $\calG\ftimes{\varphi}{\psi}\calH$ is the diffeological groupoid whose object space is the pullback of diffeological spaces $(\calG_0)\ftimes{\varphi_0}{\psi_0}(\calH_0)$, whose arrow space is the pullback of diffeological spaces $(\calG_1)\ftimes{\varphi_1}{\psi_1}(\calH_1)$, and all structure maps are restrictions of those of  the product $\calG\times\calH$.
			
			\item\labell{i:pullback gpd by f}Let $\calG$ be a diffeological groupoid and $f\colon X\to\calG_0$ a smooth map.  The \textbf{pullback of $\calG$ by $f$} is the diffeological groupoid $f^*\calG:=X^2\!\ftimes{f^2}{(\src,\trg)}\calG_1\toto X$ whose source and target maps are the first and second projections, resp.  The unit, multiplication, and inversion maps are those induced by $\calG$.  The pullback groupoid comes equipped with a smooth functor $\widehat{f}\colon f^*\calG\to\calG$ given by $(x_0,x_1,g)\mapsto g$, which is equal to $f$ on objects.
			
			\item\labell{i:relation gpd}Let $X$ be a diffeological space and $\sim$ an equivalence relation on $X$ with quotient map $\pi$.  The \textbf{relation groupoid of $\sim$} is the fibre product groupoid $X\ftimes{\pi}{\pi}X\toto X$ with source and target the projection maps $\pr_1$ and $\pr_2$ to $X$.
			
			In the extreme case that $\sim$ is equality,  then the corresponding relation groupoid has an arrow space identified with $X$ itself, with source and target the identity maps.  This is the \textbf{trivial groupoid} of $X$, denoted by $X$.  On the other extreme, if $x\sim y$ for all $x,y\in X$, then the corresponding relation groupoid is the \textbf{pair groupoid}, denoted $X^2\toto X$, whose source and target maps are the first and second projection maps, resp.  Any diffeological groupoid $\calG$ has a natural smooth functor $\chi_\calG$ to the pair groupoid $\calG_0^2\toto\calG_0$, called the \textbf{characteristic functor}, equal to the identity on objects and which sends arrows $g$ to $(\src_\calG(g),\trg_\calG(g))$.
			
			\item\labell{i:action groupoid}Let $X$ be a diffeological space and $\calG$ a diffeological groupoid acting (on the left) on $X$ with anchor map $a$ (see \cite[Definition 4.1]{vdS:morita} for a definition of a diffeological groupoid action).  The \textbf{action groupoid} of the action is the diffeological groupoid $\calG\ltimes X:=(\calG_1\ftimes{\src_\calG}{a}X\toto X)$ with source the projection $\src_{G\ltimes X}(g,x)=x$ and target the action map $\trg_{G\ltimes X}(g,x)=g\cdot x$.  A similar definition holds for right actions.
			
			\item\labell{i:neb gpd}Given a generating family $\calF$ of a diffeological space $X$, the \textbf{nebulaic groupoid of $\calF$} is the diffeological groupoid $\calN(\calF):=\ev_\calF^*X$ whose arrow space is identified with $\Neb(\calF)\ftimes{\ev_\calF}{\ev_\calF}\Neb(\calF)$.   This groupoid appears in \cite{KWW} and is used to define \c{C}ech cohomology of diffeological sheaves; it is similar to the structure groupoid appearing in \cite{iglesias-prato}. 
			
			\item\labell{i:kernel gpd}Given a smooth functor $\varphi\colon\calG\to\calH$, the \textbf{kernel groupoid of $\varphi$}, denoted $\ker(\varphi)$, is the diffeological groupoid whose arrow space is the preimage of the units of $\calH$ $$\ker(\varphi)_1=\{g\in\calG\mid\varphi(g)=\unit_{\src_\calH(\varphi(g))}\}.$$  There is a natural inclusion functor from $\ker(\varphi)$ into $\calG$ that is the identity on objects.
			
			\item\labell{i:inertia gpd}Given a diffeological groupoid $\calG$, the kernel $\ker(\chi_\calG)$ has arrow space $$\ker(\chi_\calG)_1=\{k\in\calG_1\mid \src_\calG(k)=\trg_\calG(k)\}$$ and admits a left action of $\calG$ with anchor $\trg_\calG\circ i$, where $i\colon\ker(\chi_\calG)\to\calG$ is the inclusion functor, and the action is given by conjugation $g\cdot k:= gkg^{-1}$.  The corresponding action groupoid $\calI_\calG:=\calG\ltimes\ker(\chi_\calG)_1$ is the \textbf{inertia groupoid of $\calG$}.\qedhere
		\end{enumerate}
	\end{examples}

%%%%%%%%%%%%%%%%%%%%%%%%%%%%%%%%%%%%%%%%%%%%%%%%
\section{The $2$-Site Structure on $\DGpoid$}\labell{s:2-site}
%%%%%%%%%%%%%%%%%%%%%%%%%%%%%%%%%%%%%%%%%%%%%%%%

	Since groupoids are categories, it makes sense to talk about a functor between them that is (part of) an equivalence of categories.  This is a pair of functors $F\colon \calC\to\calD$ and $G\colon\calD\to\calC$ and a pair of natural isomorphisms connecting the compositions $F\circ G$ and $G\circ F$ to the identity functors $\id_\calD$ and $\id_\calC$.  This is often called an ``adjoint equivalence''\footnotemark[1]. A weaker (but equivalent) definition that is often used is a functor that is essentially surjective and fully faithful.  Since the categories we are concerned with carry more structure, we require more from these ``weak equivalences''.  In particular, for diffeological groupoids we need smooth versions of essential surjectivity and fully faithfulness, which we will define below.  The drawback is that given such a weak equivalence $F$, we may not be able to construct a ``weak inverse'' $G$ as above with the required structure.  The solution to this is to ``formally invert'' weak equivalences, enlarging our $2$-category into a bicategory that contain formal inverses of weak equivalences.
	
	There are several common recipes for this bicategory, including that of Pronk and Scull \cite{pronk,pronk-scull:bicateg}, Roberts \cite{roberts:internal,roberts:anafunctors}, the bibundle setup \cite{HS,MM:sheaves,vdS:morita}, and stacks \cite{lerman:orbifolds,villatoro:phd}; we will focus on using the anafunctor recipe of Roberts, and later connect this to the bibundle recipe already worked through for diffeological groupoids by van der Schaaf in \cite{vdS:thesis,vdS:morita}.  
\footnotetext[1]{Some authors also require the so-called Triangle Identities as well \cite[Equation IV.1.9]{maclane}.}

	\begin{definition}[Weak Equivalence]\labell{d:weak equiv}
		Given diffeological groupoids $\calG$ and $\calH$, a smooth functor $\varphi\colon\calG\to\calH$ is 
		\begin{enumerate}
			\item\labell{i:ess surj} \textbf{smoothly essentially surjective} if the following map is a subduction:
				$$\Psi_\varphi\colon(\calG_0)\ftimes{\varphi}{\trg}\calH_1\longrightarrow\calH_0\colon (x,h)\longmapsto \src_\calH(h),$$			
			\item\labell{i:fully faithful} \textbf{smoothly fully faithful} if the following map is a diffeomorphism:
				$$\Phi_\varphi\colon \calG_1\longrightarrow(\calG_0^2)\ftimes{\varphi^2}{(\src,\trg)}\calH_1\colon g\longmapsto(\src_\calG(g),\trg_\calG(g),\varphi(g)),$$
			\item\labell{i:weak equivalence} a \textbf{weak equivalence} if it is both smoothly essentially surjective and smoothly fully faithful.
			\item\labell{i:subd weak equiv} a \textbf{subductive weak equivalence} if it is a weak equivalence and $\varphi_0$ is a subduction.
		\end{enumerate}
		We will denote a weak equivalence by $\varphi\colon\calG\we\calH$ and a subductive weak equivalence by $\varphi\colon\calG\subdwe\calH$.  Denote the class of all weak equivalences of $\DGpoid$ by $\calW$ and the class of all subductive weak equivalences of $\DGpoid$ by $\calJ$.
	\end{definition}
	
	It turns out requiring the ``smoothness'' of our weak equivalence now separates this concept from the adjoint equivalence defined above, as the following example illustrates.
	
	\begin{example}[Principal Bundles]\labell{x:princ bdle}
		Let $G$ be a Lie group and $\pi\colon P\to B$ a principal $G$-bundle over a manifold $B$. Consider the Lie groupoid given by the action groupoid $G\ltimes P$, and the smooth functor $\Pi\colon G\ltimes P\to B$ sending $(g,p)$ to $\pi(p)$.
		
		Let $q$ be a plot of $B$.  After shrinking $U_q$, there exist an open subset $V$ containing $q(U_q)$ and a smooth section $\sigma\colon V\to P|_V$.  The plot $r=(\sigma\circ q,q)$ is a lift of $q$ from $B$ to $P\ftimes{\pi}{\trg}B$.  It follows that $\Psi_\Pi$ is a subduction.
		
		Since the $G$-action is free, $\Phi_\Pi$ is injective.  Since the division map $P\ftimes{\pi}{\pi}P\to G$ sending $(p_1,p_2)$ to the unique $g$ such that $p_2=g\cdot p_1$ is smooth, it follows that $\Phi_\Pi$ is subductive. Therefore, $\Pi$ is smoothly fully faithful, and hence a weak equivalence.  Since $\pi$ is a subduction, $\Pi$ is a subductive weak equivalence (and checking that $\Psi_\Pi$ is subductive was in fact unnecessary; see \cref{r:weak equiv} below).
		
		This example is useful for illustrating the difference between (smooth) weak equivalency and adjoint equivalency: the reader may check that $\Pi$ is an adjoint equivalence if and only if $\pi\colon P\to B$ is a trivial bundle (the required functor from $B$ to $G\ltimes P$ induces a global section of $\pi$).
	\end{example}
	
	The notation $\Psi$ and $\Phi$ in \cref{d:weak equiv} were inspired by the words ``surjective'' and ``fully faithful'', respectively, which may help the reader recall what they are in the sequel.
	
	\begin{remark}\labell{r:weak equiv}
		Any smoothly fully faithful functor $\varphi\colon\calG\to\calH$ such that $\varphi_0$ is a subduction is smoothly essentially surjective, and hence $\varphi$ is a subductive weak equivalence.
	\end{remark}
	
	\begin{remark}\labell{r:morita fibration}
		In \cite{dHF}, the authors define a ``Morita fibration'' $\varphi\colon\calG\to\calH$ between Lie groupoids to be a (Lie) weak equivalence that is also a ``fibration''; that is, both $\varphi_0$ and the map $\operatorname{fibr}\colon\calG_1\to\calH_1\times_{\calH_0}\calG_0$ sending $g$ to $(\varphi(g),\src(g))$ are surjective submersions.  In fact, that $\fibr$ is surjective submersive is automatic given a weak equivalence that is surjective submersive on objects.  Thus in the Lie setting, ``surjective submersive weak equivalences'' (weak equivalences that are surjective submersive on objects) are exactly Morita fibrations.  Subductive weak equivalences between diffeological groupoids are a generalisation of this: given a subductive weak equivalence, the corresonding map $\fibr$ is subductive.  Moreover, if $\varphi$ is locally subductive, then so is $\fibr$.
		
		The authors in \cite{dHF} also indicate that a smooth functor $\varphi\colon\calG\to\calH$ between Lie groupoids is a Morita fibration if and only if its kernel is isomorphic to the Lie groupoid $\calG_0\ftimes{\varphi}{\varphi}\calG_0\toto\calG_0$.  A similar statement also holds for subductive weak equivalences: $\varphi\colon\calG\to\calH$ is a subductive weak equivalence if and only if $\ker(\varphi)$ is isomorphic to $(\calG_0)\ftimes{\varphi}{\varphi}(\calG_0)$ as diffeological groupoids. 
	\end{remark}
	
	Our goal is to ``formally invert'' the elements of $\calW$ via a ``localisation'' of $\DGpoid$ at $\calW$ (to be defined later).  Roberts in \cite{roberts:anafunctors} shows that the anafunctor recipe for localisation can be achieved provided that the $2$-category $\DGpoid$ can be equipped with a $2$-site structure. We need the following terminology \cite[Definitions 2.2, 2.9, 2.12]{roberts:anafunctors}:
	
	\begin{definition}[Singleton Strict Pretopology]\labell{d:pretop}
		Let $\calB$ be a bicategory.
		\begin{enumerate}
			\item\labell{i:ff} A $1$-cell $f\colon b\to c$ of $\calB$ is \textbf{representably fully faithful} if for any $1$-cells $g,h\colon a\to b$ and $2$-cell $A\colon f\circ g\Rightarrow f\circ h$, there is a unique $2$-cell $A'\colon g\Rightarrow h$ such that $A = fA'$. 
			\item\labell{i:co-ff} A $1$-cell $f\colon a\to b$ is \textbf{co-fully faithful} if for any $1$-cells $g,h\colon b\to c$ and $2$-cell $A\colon g\circ f\Rightarrow h\circ f$, there is a unique $2$-cell $A'\colon g\Rightarrow h$ such that $A=A'f$.
			\item\labell{i:pretop} A class $C$ of $1$-cells in $\calB$ is a \textbf{singleton strict pretopology} if it contains all of the identity $1$-cells; is closed under composition; and for any $1$-cells $f\colon a\to b$ of $\calB$ and $g\colon c\to b$ in $C$, the pullback
	$$\xymatrix{
		a\ftimes{f}{g} c \ar[r] \ar[d]_{h} & c \ar[d]^{g} \\
		a \ar[r]_{f} & b \\
	}$$
exists with $h\in C$.
			\item\labell{i:bi-ff} A singleton strict pretopology $C$ of $\calB$ is \textbf{bi-fully faithful} if every $f\in C$ is both representably fully faithful and co-fully faithful.
		\end{enumerate}
		We call $\calB$ a \textbf{$2$-site} after equipping it with a bi-fully faithful singleton strict pretopology, denoted $(\calB,C)$.
	\end{definition}
	
	The main work involved in proving that $(\DGpoid,\calJ)$ is a $2$-site is contained in the following lemma.
	
	\begin{lemma}\labell{l:we properties}
		\noindent
		\begin{enumerate}
			\item\labell{i:FF to fully faithful} A smooth functor is representably fully faithful if and only if it is smoothly fully faithful.
			\item\labell{i:J co-FF} A subductive weak equivalence is co-fully faithful.
			\item\labell{i:we comp} Given smooth functors $\varphi\colon\calG\to\calH$ and $\psi\colon\calH\to\calK$, if any two of $\varphi$, $\psi$, and $\psi\circ\varphi$ are weak equivalences, then so is the third.
			\item\labell{i:subdwe pullback} Given the following pullback diagram, 
				$$\xymatrix{
					\calG\ftimes{\varphi}{\psi}\calH \ar[d]_{\pr_1} \ar[r]^{\quad\pr_2} & \calH \ar[d]^{\psi} \\
					\calG \ar[r]_{\varphi} & \calK \\
				}$$
				if $\varphi$ is a subductive weak equivalence, then so is $\pr_2$.
		\end{enumerate}
	\end{lemma}
	
	\begin{proof}
		Suppose $\varphi\colon\calG\we\calH$ is smoothly fully faithful; $\psi,\chi\colon\calK\to\calG$ are smooth functors; and $S\colon\varphi\circ\psi\Rightarrow\varphi\circ\chi$ is a smooth natural transformation.  Define $S'\colon\psi\Rightarrow\chi$ by $$S'_z:=\Phi_\varphi^{-1}(\psi(z),\chi(z),S_z).$$  Since $\varphi$ is smoothly fully faithful, $S'$ is well-defined and smooth, is a natural transformation, and is the unique natural transformation satisfying $S=\varphi S'$.

		Conversely, suppose that for any smooth functors $\psi,\chi\colon\calK\to\calG$ and smooth natural transformation $S\colon\varphi\circ\psi\Rightarrow\varphi\circ\chi$, there exists a unique natural transformation $S'\colon\psi\Rightarrow\chi$ such that $S=\varphi S'$.  Fix a point $(x,x',h)\in(\calG_0^2)\ftimes{\varphi^2}{(\src,\trg)}\calH_1$, and set $\calK$ to be the trivial groupoid of a point, $\psi_0$ to have image $x$, $\chi_0$ to have image $x'$, and $S$ to have image $h$.  There is a unique smooth $S'\colon\psi\Rightarrow\chi$ such that $\varphi S'=S$, from which it follows that there is a unique $g\in\calG_1$ with source $x$, target $x'$, and such that $\varphi(g)=h$.  So $\Phi_\varphi$ is bijective.

		Fix a plot $p=(p_1,p_2,p_3)\colon u\mapsto (x_u,x'_u,h_u)$ of $(\calG_0^2)\ftimes{\varphi^2}{(\src,\trg)}\calH_1$, set $\calK$ to be the trivial groupoid of $U_p$, $\psi_0=p_1$, $\chi_0=p_2$, and $S=p_3$.  There is a unique smooth $S'\colon\psi\Rightarrow\chi$ such that $\varphi S'=S$, from which it follows that $\Phi_\varphi\circ S'=p$.  It follows that $\Phi_\varphi$ is a subduction.  Since any injective subduction is a diffeomorphism, $\varphi$ is smoothly fully faithful.  This proves \cref{i:FF to fully faithful}.
		
		Let $\varphi\colon\calG\subdwe\calH$ be a subductive weak equivalence; $\psi,\chi\colon\calH\to\calK$ smooth functors; and $S\colon\psi\circ\varphi\Rightarrow\chi\circ\varphi$ a smooth natural transformation.  Define $S'\colon\calH_0\to\calK_1$ by $S'_y=S_x$ where $x\in\varphi_0^{-1}(y)$.  To show that this is well-defined, suppose $\varphi(x_1)=\varphi(x_2)$.  Since $\varphi$ is a weak equivalence, there exists $g=\Phi_\varphi^{-1}(x_1,x_2,\unit_{\varphi(x_1)})$ from $x_1$ to $x_2$, inducing the commutative diagram
		$$\xymatrix{
			\psi\circ\varphi(x_1) \ar[d]_{\psi\circ\varphi(g)} \ar[r]^{S_{x_1}} & \chi\circ\varphi(x_1) \ar[d]^{\chi\circ\varphi(g)} \\
			\psi\circ\varphi(x_2) \ar[r]_{S_{x_2}} & \chi\circ\varphi(x_2). \\
		}$$
However, since $\varphi(g)=\unit_{\varphi(x_1)}$, we have $S_{x_1}=S_{x_2}$.  This shows that $S'$ is well-defined on the image of $\varphi$, which is $\calH_0$ since $\varphi$ is subductive.  By construction, $S'$ is the unique natural transformation satisfying $S=S'\varphi$.

		To show that $S'$ is smooth, let $p$ be a plot of $\calH_0$.  Since $\varphi$ is subductive, after shrinking $U_p$, there exists a lift $q$ of $p$ to $\calG_0$.  Then $S'\circ p=S\circ q$, the latter of which is a plot of $\calK_1$.  This proves \cref{i:co-ff}.
		
		\cref{i:we comp} is known as the ``3 for 2 property'' and follows from a similar property for fibre products.  See \cite[Lemma 8.1]{pronk-scull:translation-err} for a proof using Lie groupoids; there, Lie groupoids can be replaced with diffeological groupoids, and surjective submersions with subductions.
		
		By \cref{i:pullback subd} of \cref{x:ind subd}, to prove \cref{i:subdwe pullback}, it remains to show that $\pr_2$ is smoothly fully faithful.  Since $\varphi$ is smoothly fully faithful, given $(g,h)\in(\calG\ftimes{\varphi}{\psi}\calH)_1$ we have $g=\Phi^{-1}_{\varphi}(\src_\calG(g),\trg_\calG(g),\psi(h))$.  It follows that $\Phi_{\pr_2}$ is injective and subductive, hence a diffeomorphism.
	\end{proof}

	The following proposition now follows from \cref{i:comp ind subd} of \cref{x:ind subd} and \cref{l:we properties}.
	
	\begin{proposition}[$(\Diffeol,\calJ)$ is a $2$-Site]\labell{p:dgpoid 2-site}
		The class $\calJ$ is a bi-fully faithful singleton strict pretopology on $\DGpoid$, making $(\DGpoid,\calJ)$ a $2$-site.
	\end{proposition}

%%%%%%%%%%%%%%%%%%%%%%%%%%%%%%%%%%%%%%%%%%%%%%%%
\section{The Anafunctor Bicategory}\labell{s:anafunctors}
%%%%%%%%%%%%%%%%%%%%%%%%%%%%%%%%%%%%%%%%%%%%%%%%		

	Following \cite{roberts:anafunctors}, we now construct a bicategory of diffeological groupoids whose $1$-cells are so-called ``anafunctors'', also called $\calJ$-fractions in \cite{roberts:anafunctors}.  

	\begin{definition}[Anafunctor]\labell{d:anafunctor}
		An \textbf{anafunctor} is a pair of smooth functors $\calG\ot{\varphi}\calK\to\calH$ in which the left functor $\varphi$ is a subductive weak equivalence.  The \textbf{identity anafunctor} of a diffeological groupoid $\calG$ is the anafunctor $\calG\overset{=}{\longleftarrow}\calG\overset{=}{\longrightarrow}\calG$.
	\end{definition}

The composition of two anafunctors uses a strict pullback groupoid to define it; see \cref{i:pullback gpd} of \cref{x:diffeol gpd}.

	\begin{definition}[Composition of Anafunctors]\labell{d:anaf comp}
		Let $\calG\ot{\varphi}\calL\uto{\psi}\calH$ and $\calH\ot{\chi}\calM\uto{\omega}\calK$ be anafunctors.  Define their \textbf{composition} to be  the anafunctor $\calG\ot{\varphi\circ\pr_1}\calL\ftimes{\psi}{\chi}\calM\uto{\omega\circ\pr_2}\calK$.  
	\end{definition}

The $2$-cells of the bicategory are certain natural transformations, a feature which makes this bicategory friendlier than equivalent bicategories of groupoids.

	\begin{definition}[Transformation]\labell{d:transformation}
		Given anafunctors $\calG\ot{\varphi}\calK\uto{\psi}\calH$ and $\calG\ot{\varphi'}\calK'\uto{\psi'}\calH$, a \textbf{transformation} between them is a natural transformation
		$$\xymatrix{\ar @{} [dr] |{\rotatebox{45}{$\Rightarrow$}S}
			\calK\ftimes{\varphi}{\varphi'}\calK' \ar@{>>}[d]_{\pr_1}^{\simeq} \ar@{>>}[r]^{\pr_2}_{\simeq} & \calK' \ar[d]^{\psi'} \\
			\calK \ar[r]_{\psi} & \calH.
		}$$
		The \textbf{identity transformation} of an anafunctor $\calG\ot{\varphi}\calK\uto{\psi}\calH$ is given by the natural transformation $$I_{\calG\leftarrow\calK\to\calH}\colon(\calK\ftimes{\varphi}{\varphi}\calK)_0\to\calH_1\colon(y_1,y_2)\mapsto\psi(\Phi_\varphi^{-1}(y_1,y_2,\unit_{\varphi(y_1)})).\eqno\qedhere$$
	\end{definition}

	Vertical and horizontal composition of transformations, along with associators and unitors, we leave undefined as we do not use them here, and instead refer to \cite{roberts:anafunctors} for their definitions. We now can apply \cite[Proposition 3.20]{roberts:anafunctors} to our setting; since $(\DGpoid,\calJ)$ is a $2$-site by \cref{p:dgpoid 2-site}, we have:

	\begin{proposition}[Bicategory of Diffeological Groupoids]\labell{p:bicateg diffeol gpds}
		Diffeological groupoids form a bicategory $\DGpoid_\ana$ with anafunctors as $1$-cells and transformations of anafunctors as $2$-cells.
	\end{proposition}
	
	What is special about $\DGpoid_\ana$ is that it comes with an inclusion of $\DGpoid$ into it, and also provides inverses for weak equivalences after passing through this inclusion.  More specifically, recall that an \textbf{equivalence} in a bicategory $\calB$ is a $1$-cell $F\colon a\to b$ that has a \textbf{quasi-inverse} $\bar{F}\colon b\to a$ for which the composition $F\circ\bar{F}$ admits an invertible $2$-cell $\eps_F$ to $\id_b$ and $\bar{F}\circ F$ admits an invertible $2$-cell $\eta_F$ to $\id_a$.  A \textbf{localisation} of a bicategory $\calB$ with respect to a class of $1$-cells $C$ is a bicategory $\widetilde{\calB}$ and a pseudofunctor $L\colon \calB\to\widetilde{\calB}$ such that all elements of $L(C)$ are equivalences, and $L$ is universal in the sense that precomposition with $L$ induces an equivalence of bicategories $$L^*\colon\Bicat(\widetilde{\calB},\calA)\to\Bicat_C(\calB,\calA),$$ where $\Bicat_C$ is the full sub-bicategory on the pseudofunctors sending $C$ to equivalences in $\calA$.  In our case, $C=\calW$ and a choice of $L$ is given by spanisation:
	
	\begin{definition}[Spanisation]\labell{d:spanisation}
		Given a smooth functor $\varphi\colon\calG\to\calH$, the \textbf{spanisation} of $\varphi$ is the anafunctor $\calG\overset{=}{\longleftarrow}\calG\uto{\varphi}\calH$.  Given a smooth natural transformation
	$$\xymatrix{
		\calG \ar@/^1pc/[rr]^{\varphi} \ar@/_1pc/[rr]_{\psi} & \Downarrow S & \calH \\
	}$$
	the \textbf{spanisation of $S$} is the transformation
	$$\begin{gathered}[b]\xymatrix{
		\calG\ftimes{\id}{\id}\calG \ar[r]_{\quad\cong} & \ar@{}[dr] |{\rotatebox{45}{$\Leftarrow$}S}\calG \ar[r]^{=} \ar[d]_{=} & \calG \ar[d]^{\varphi} \\
		 & \calG \ar[r]_{\psi} & \calH\\
	}\\[-\dp\strutbox]\end{gathered}\eqno\qedhere$$
	\end{definition}
	
	\begin{theorem}[Spanisation is a Localisation of $\DGpoid_\ana$]\labell{t:spanisation}
		The pseudofunctor $\mathfrak{S}\colon\DGpoid\to\DGpoid_\ana$ sending smooth functors and smooth natural transformations to their spanisations is a localisation of $\DGpoid$ at $\calW$.
	\end{theorem}
	
	This theorem is essentially already proven by \cite[Theorem 3.24]{roberts:anafunctors} given \cref{p:dgpoid 2-site,p:bicateg diffeol  gpds}; the only part missing is that we need to check that weak equivalences are so-called $\calJ$-locally split (see \cite[Definition 3.22]{roberts:anafunctors}).
	
	\begin{definition}[$\calJ$-Locally Split]\labell{d:loc split}
		A smooth functor $\varphi\colon\calG\to\calH$ is \textbf{$\calJ$-locally split} if there is a subductive weak equivalence $\psi\colon\calK\onto\calH$, a smooth functor $\chi\colon\calK\to\calG$ and a smooth natural transformation
		$$\begin{gathered}[b]\xymatrix{\ar@{}[dr] |{\quad\Downarrow S}
			 & \calG \ar[d]^{\varphi} \\
			\calK \ar@/^/[ur]^{\chi} \ar@{>>}[r]_{\psi}^{\cong} & \calH.\\
		}\\[-\dp\strutbox]\end{gathered}\eqno\qedhere$$
	\end{definition}
	
	The relationship between $\calJ$-locally split functors and weak equivalences is given by the following lemma.
	
	\begin{lemma}\labell{l:loc split}
		A smooth functor $\varphi\colon\calG\to\calH$ is a weak equivalence if and only if it is smoothly fully faithful and $\calJ$-locally split.
	\end{lemma}
	
	In order to prove this lemma, we introduce the weak pullback and one of its important properties. 
	
	\begin{definition}[Weak Pullback]\labell{d:weak pullback}
		Let $\varphi\colon\calG\to\calK$ and $\psi\colon\calH\to\calK$ be smooth functors.  Define the \textbf{weak pullback} of $\varphi$ and $\psi$ to be the diffeological groupoid $\calG\wtimes{\varphi}{\psi}\calH$ in which
		$$(\calG\wtimes{\varphi}{\psi}\calH)_0:=(\calG_0)\ftimes{\varphi}{\src}(\calK_1)\ftimes{\trg}{\psi}(\calH_0),$$
		the space of triples $(x,k,y)\in\calG_0\times\calK_1\times\calH_0$ such that $k$ is an arrow from $\varphi(x)$ to $\psi(y)$; and
		$$(\calG\wtimes{\varphi}{\psi}\calH)_1:=(\calG_1)\ftimes{\src\circ\varphi}{\src}(\calK_1)\ftimes{\trg}{\src\circ\psi}(\calH_1),$$
		the space of triples $(g,k,h)\in\calG_1\times\calK_1\times\calH_1$ with source and target maps
		$$\src(g,k,h)=(\src_\calG(g),k,\src_\calH(h)) \quad \text{and} \quad \trg(g,k,h)=(\trg_\calG(g),\psi(h)\circ k\circ\varphi(g)^{-1},\trg_\calH(h)).$$
		Thus, given triples $(x_1,k_1,y_1),(x_2,k_2,y_2)\in(\calG\wtimes{\varphi}{\psi}\calH)_0$, an arrow $(g,k_1,h)$ from $(x_1,k_1,y_1)$ to $(x_2,k_2,y_2)$ induces a commutative diagram
		$$\xymatrix{
			\varphi(x_1) \ar[r]^{k_1} \ar[d]_{\varphi(g)} & \psi(y_1) \ar[d]^{\psi(h)} \\
			\varphi(x_2) \ar[r]_{k_2} & \psi(y_2). \\
		}$$
		
		The other structure maps are given as follows.  The unit map is given by $$\unit(x,k,y)=(\unit_\calG(x),k,\unit_\calH(y)),$$ multiplication is given by
		$$\mult((g_2,k_2,h_2),(g_1,k_1,h_1))=(g_2g_1,k_1,h_2h_1)$$
		for composable triples $(g_1,k_1,h_1)$ and $(g_2,k_2,h_2)$, and inversion is given by $$\inv(g,k,h)=(g^{-1},h\cdot k\cdot g^{-1},h^{-1}).\eqno\qedhere$$
	\end{definition}
	
	\begin{remark}\labell{r:weak pullback}
		The diagram induced by the weak pullback is $2$-commutative, with the $2$-cell the natural transformation $\PR_2\colon\varphi\Rightarrow\psi$ given by $(\calG\wtimes{\varphi}{\psi}\calH)_0\to\calK_1\colon(x,k,y)\mapsto k$:
		$$\begin{gathered}[b]\xymatrix{\ar@{}[dr] |{\quad\rotatebox{45}{$\Rightarrow$}\Pr_2}
			\calG\wtimes{\varphi}{\psi}\calH \ar[r]^{\pr_3} \ar[d]_{\pr_1} & \calH \ar[d]^{\psi} \\
			\calG \ar[r]_{\varphi} & \calK \\
		}\\[-\dp\strutbox]\end{gathered}\eqno\qedhere$$
	\end{remark}
	
	\begin{lemma}\labell{l:weak pullback equiv}
		Given smooth functors $\varphi\colon\calG\to\calK$ and $\psi\colon\calH\to\calK$, if $\psi$ is a weak equivalence, then $\pr_1\colon\calG\wtimes{\varphi}{\psi}\calH\to\calH$ is a subductive weak equivalence.
	\end{lemma}
	
	\begin{proof}
		Suppose $\psi$ is a weak equivalence.  Fix a plot $p\colon u\mapsto x_u$ of $\calG_0$.  By smooth essential surjectivity, after shrinking $U_p$, there is a lift $u\mapsto (y_u,k_u),$ of $\varphi\circ p$ to $(\calH_0)\ftimes{\psi}{\trg}\calK_1$.  Then $u\mapsto(x_u,k_u,y_u)$ is a lift of $p$ to $(\calG\wtimes{\varphi}{\psi}\calH)_0$, and so $(\pr_1)_0$ is a subduction.

		Since $\Phi_\psi$ is a diffeomorphism, $\Phi_{\pr_1}$ is injective.  Moreover, a plot $$p\colon u\mapsto ((x_u,k_u,y_u),(x'_u,k'_u,y'_u),g_u)$$ of $\left((\calG\wtimes{\varphi}{\psi}\calH)_0^2\right)\ftimes{\pr_1^2}{(\src,\trg)}\calG_1$ has a lift to $(\calG\wtimes{\varphi}{\psi}\calH)_1$ given by $$u\mapsto(g_u,k_u,\Phi_\psi^{-1}(y_u,y'_u,k'_u\cdot\varphi(g_u)\cdot k_u^{-1})).$$
It follows that $\Phi_{\pr_1}$ is a subduction, and hence a diffeomorphism.  This proves that $\pr_1$ is a subductive weak equivalence.
	\end{proof}
	
	We are now ready to prove \cref{l:loc split}.
	
	\begin{proof}[Proof of \cref{l:loc split}]
		Suppose $\varphi$ is a weak equivalence.  Choose $\calK:=\calG\wtimes{\varphi}{\id_\calH}\calH$, $\psi=\pr_3$, $\chi=\pr_1$, and $S=\PR_2$.  Then $\psi$ is a subductive weak equivalence by \cref{l:weak pullback equiv}.

		Conversely, suppose $\varphi$ is smoothly fully faithful and that there exist a subductive weak equivalence $\psi\colon\calK\to\calH$, a functor $\chi\colon\calK\to\calG$, and a natural transformation $S\colon\varphi\circ\chi\Rightarrow\psi$.  Let $p$ be a plot of $\calH_0$.  Since $\psi_0$ is a subduction, after shrinking $U_p$, there is a lift $q$ of $p$ to $\calK_0$.  Then $(\chi\circ q,(S\circ q)^{-1})$ is the desired lift of $p$ to $\calG_0\ftimes{\varphi}{\trg_\calH}\calH_1$, implying that $\Psi_\varphi$ is a subduction.
	\end{proof}
	
	With this lemma, we have proven \cref{t:spanisation}.  Another use of the weak pullback is that it allows us to define a quasi-inverse of the spanisation of a weak equivalence quite easily.
	
	\begin{proposition}[Quasi-Inverses to Weak Equivalences]\labell{p:quasi-inverse}
		Given a weak equivalence $\varphi\colon\calG\we\calH$, a quasi-inverse to its spanisation $\mathfrak{S}(\varphi)$ is the anafunctor $\calH\ot{\pr_1}\calH\wtimes{\id_\calH}{\varphi}\calG\overset{\simeq}{\underset{\raisebox{1mm}{$\scriptstyle\pr_3$}}{\onto}}\calG$.
	\end{proposition}
	
	\begin{proof}
		It is straightforward to check that the natural transformation $\PR_2$ induced by the weak pullback (see \cref{r:weak pullback}) induces transformations between each of the two compositions of these anafunctors to the required identity anafunctors.
	\end{proof}
	
	\begin{remark}\labell{r:quasi-inverse}
		One may notice that the right map $\pr_3$ to $\calG$ of the quasi-inverse in \cref{p:quasi-inverse} is subductive, even though $\varphi$ is not.  One may ask whether there is a $2$-cell from $\mathfrak{S}(\varphi) = \left(\calG\ot{\id_\calG}\calG\underset{\varphi}{\we}\calH\right)$ to an anafunctor whose right arrow is subductive.  The answer is yes: the anafunctor $\calG\ot{\pr_1}\calG\wtimes{\varphi}{\id_\calH}\calH\overset{\simeq}{\underset{\raisebox{1mm}{$\scriptstyle\pr_3$}}{\onto}}\calH$.  This is called the \textbf{anafunctisation} of $\varphi$ in \cite{duli}.  Note that it is the mirror image of the quasi-inverse from \cref{p:quasi-inverse}.
	\end{remark}
	
	We end this section with a definition of Morita equivalence in terms of anafunctors.
	
	\begin{definition}[Morita Equivalence]\labell{d:morita}
		Two diffeological groupoids are \textbf{Morita equivalent} if there is an anafunctor (called a \textbf{Morita equivalence}) between them in which both arrows are weak equivalences.
	\end{definition}
	
	\begin{remark}\labell{r:morita}
		Given a Morita equivalence, one can arrange for both arrows to be subductive weak equivalences, generalising \cref{r:quasi-inverse}.  The fact that any Morita equivalence admits a $2$-cell to either an element of $\mathfrak{S}(\calW)$ or a quasi-inverse as in \cref{p:quasi-inverse} follows from the universal property of localisation.
	\end{remark}

%%%%%%%%%%%%%%%%%%%%%%%%%%%%%%%%%%%%%%%%%%%%%%%%
\section{Lie Groupoids}\labell{s:lie gpds}
%%%%%%%%%%%%%%%%%%%%%%%%%%%%%%%%%%%%%%%%%%%%%%%%	

	We now restrict our attention to Lie groupoids, with the goal of showing that the anafunctor bicategory of Lie groupoids is an essentially full sub-bicategory of $\DGpoid_\ana$.  There is subtlety here that the casual reader may initially miss: weak equivalences in the Lie setting require more than those in the diffeological setting.  This is alleviated by some perhaps surprising facts given in \cref{l:weak equiv local subd}.
	
	In the definition of Lie groupoid, it is common to relax the Hausdorff and second-countability criteria for the arrow space; indeed, such Lie groupoids occur naturally in the study of foliations as holonomy groupoids, or groupoids of germs of local diffeomorphisms.  Therefore, we do not require these conditions.  However, one \emph{could} require either one or both, and the results of this paper still hold.  The key to this is the fact that these properties pull back by weak equivalences; see \cref{p:pullback}.

	\begin{definition}[Lie Groupoid]\labell{d:lie gpd}
		A \textbf{Lie groupoid} is a diffeological groupoid $\calG$ in which $\calG_0$ and $\calG_1$ are smooth manifolds with $\src_\calG$ and $\trg_\calG$ submersions.   We do not require $\calG_1$ to be Hausdorff and second-countable.  However, if $\calG_1$ is Hausdorff (resp.\ second-countable), then we say that $\calG$ is \textbf{Hausdorff} (resp.\ \textbf{second-countable}).  Lie groupoids form a full sub-2-category $\LieGpoid$ of $\DGpoid$.
	\end{definition}
	
	\begin{example}[Pullback by Surjective Submersion]\labell{x:pullback surj subm}
		Let $\calG$ be a Lie groupoid, and let $f\colon M\to\calG_0$ be a surjective submersion from a manifold $M$.  The \textbf{pullback groupoid by $f$}, denoted $f^*\calG$, is the groupoid $(M^2)\ftimes{f^2}{(\src,\trg)}\calG_1\toto M$, whose source and target maps are the first and second projections, resp.  The unit, multiplication, and inversion maps are those induced by $\calG$.  By \cite[Subsection 1.4]{MM:sheaves}, $f^*\calG$ is a Lie groupoid if $\Psi_{\widehat{f}}$ is submersive (the domain is a manifold since $\trg_\calG$ is a surjective submersion); here, $\widehat{f}$ is the functor $(f,\pr_3)$.
		
		The surjectivity of $\Psi_{\widehat{f}}$ follows from the surjectivity of $f$.  Fix a plot $p\colon u\mapsto x_u$ of $\calG_0$ and a point $(w_0,g_0)\in M\ftimes{f}{\trg}\calG$ such that $\Psi_{\widehat{f}}(w_0,g_0)=x_0$.  After shrinking $U_p$, there is a lift $q\colon u\mapsto w_u$ against $f$ through $w_0$, and a lift $r\colon u\mapsto g_u$ against $\trg_\calG$ through $g_0$.  The plot $(q,r)$ of $M\ftimes{f}{\trg}\calG$ is a lift of $p$ against $\Psi_{\widehat{f}}$ through $(w_0,g_0)$.  Thus $\Psi_{\widehat{f}}$ is a surjective submersion, and we conclude that $f^*\calG$ is a Lie groupoid.
	\end{example}
	
	The localisation of $\LieGpoid$ at weak equivalences is a standard setting in differential geometry in which one can search for stacky invariants; see \cite{HS,pronk,MM:sheaves,lerman:orbifolds}.  As mentioned above, a weak equivalence in this setting \emph{a priori} is slightly different than the definition of a weak equivalence in $\DGpoid$.

	\begin{definition}[Lie and Surjective Submersive Weak Equivalences]\labell{d:surj subm weak equiv}
		 A weak equivalence $\varphi\colon\calG\to\calH$ between Lie groupoids is \textbf{Lie} if $\Phi_\varphi$ and $\Psi_\varphi$ are smooth maps between manifolds with $\Psi_\varphi$ a submersion.  A weak equivalence between Lie groupoids is \textbf{surjective submersive} if it is Lie and $\varphi_0$ is surjective submersive.  Denote by $\calW_\Lie$ the class of Lie weak equivalences in $\LieGpoid$, and by $\calJ_\Lie$ the class of surjective submersive weak equivalences in $\LieGpoid$.  An anafunctor $\calG\ot{\varphi}\calK\uto{\psi}\calH$ between two Lie groupoids is \textbf{Lie} if $\varphi$ is surjective submersive; in particular, $\calK$ is Lie.
	\end{definition}
	
	While we do not require Lie groupoids to have Hausdorff nor second-countable arrow spaces, if they happen to have one of these properties, then these properties pull back by weak equivalences.
	
	\begin{proposition}[Hausdorffness and Second-Countability Pull Back]\labell{p:pullback}
		Let $\calH$ be a Hausdorff (resp.\ second-countable) Lie groupoid.  Let $\varphi\colon\calG\to\calH$ be a weak equivalence.  Then $\calG$ is Hausdorff (resp.\ second-countable).
	\end{proposition}

	\begin{proof}
		Suppose $\calH$ is Hausdorff.  Fix $g_1,g_2\in\calG$.  If $\src_\calG(g_1)\neq \src_\calG(g_2)$ or $\trg_\calG(g_1)\neq\trg_\calG(g_2)$, then there exist disjoint open neighbourhoods of $g_1$ and $g_2$ since $\calG_0$ is Hausdorff.  Suppose $\src_\calG(g_1)=\src_\calG(g_2)$ and $\trg_\calG(g_1)=\trg_\calG(g_2)$.  By smooth fully faithfulness, $\varphi(g_1)\neq\varphi(g_2)$.  Since $\calH_1$ is Hausdorff, the preimages of disjoint open neighbourhoods of $\varphi(g_1)$ and $\varphi(g_2)$ are open and disjoint, separating $g_1$ and $g_2$.  Thus, $\calG$ is Hausdorff.
		
		Suppose $\calH$ is second-countable.  Since $\calG_0$ is second-countable, the product $\calG_0^2\times\calH_1$ is also second-countable, and hence so is the fibred product $\calG_0^2\ftimes{\varphi^2}{(\src,\trg)}\calH_1$.  By smooth fully faithfulness, $\calG$ is second-countable.
	\end{proof}

	It turns out that the slight differences between the definitions of (subductive) weak equivalence and those in \cref{d:surj subm weak equiv} are illusory.  We require a lemma.
	
	\begin{lemma}\labell{l:weak equiv local subd}
		Let $\varphi\colon\calG\to\calH$ be a weak equivalence between diffeological groupoids.
		\begin{enumerate}
			\item\labell{i:local subd1} The maps $\src_\calG$ and $\trg_\calG$ are local subductions if and only if $\Psi_\varphi$ is a local subduction.
			\item\labell{i:local subd2} If $\varphi$ is a subductive weak equivalence, then $\src_\calG$ and $\trg_\calG$ are locally subductive if and only if $\src_\calH$, $\trg_\calH$, and $\varphi_0$ are. 
		\end{enumerate}
	\end{lemma}

	\begin{proof}
		Suppose $\src_\calG$ and $\trg_\calG$ are locally subductive.  Let $p\colon u\mapsto y_u$ be a plot of $\calH_0$ and fix $(x_0,h_0)\in\calG_0\ftimes{\varphi}{\trg}\calH_1$ such that $\Psi_\varphi(x_0,h_0)=y_0$.  Since $\Psi_\varphi$ is subductive, after shrinking $U_p$, there exists a lift $u\mapsto(x'_u,h'_u)$ of $p$ to $(\calG_0)\ftimes{\varphi}{\trg}\calH_1$.  Let $g_0:=\Phi_\varphi^{-1}(x_0,x'_0,h'_0h_0^{-1})$.  Since $\trg_\calG$ is locally subductive, after shrinking $U_p$ again, there is a lift $u\mapsto g_u$ of $x'_u$ against $\trg_\calG$ through $g_0$.  Then $$r\colon U_p\to(\calG_0)\ftimes{\varphi}{\trg}\calH_1\colon u\mapsto\left(\src_\calG(g_u),\varphi(g_u)^{-1}h'_u\right)$$
is a lift of $p$ through $(x_0,h_0)$.  Thus $\Psi_\varphi$ is a local subduction.

		Conversely, suppose $\Psi_\varphi$ is locally subductive.  Fix a plot $p\colon u\mapsto x_u$ of $\calG_0$ and $g_0\in\calG_1$ such that $\src_\calG(g_0)=x_0$.   Let $(x'_0,h'_0):=(\trg_\calG(g_0),\varphi(g_0))\in(\calG_0)\ftimes{\varphi}{\trg)}\calH_1$.  Then $\Psi_\varphi(x'_0,h'_0)=\varphi(x_0)$, and after shrinking $U_p$, there is a lift $(x'_u,h'_u)$ of $\varphi\circ p$ to $(\calG_0)\ftimes{\varphi}{\trg}\calH_1$ through $(x'_0,h'_0)$.  Thus there is a lift $u\mapsto\Phi_\varphi^{-1}(x_u,x'_u,h'_u)$ of $p$ against $\src_\calG$.  But $\Phi^{-1}_\varphi(x_0,x'_0,h'_0)=g_0$.  Thus $\src_\calG$ is a local subduction (and thus so is $\trg_\calG$ since $\inv_\calG$ is a diffeomorphism).  This proves \cref{i:local subd1}.
		
		Now suppose $\varphi$ is a subductive weak equivalence, and that $\src_\calG$ and $\trg_\calG$ are locally subductive. Let $p\colon u\mapsto y_u$ be a plot of $\calH_0$ and let $h_0\in\calH_1$ such that $\src_\calH(h_0)=y_0$.  Since $\varphi$ is subductive, after shrinking $U_p$, there is a lift $x_u$ of $y_u$ to $\calG_0$; there is also some $x'_0\in\calG_0$ such that $\varphi(x'_0)=\trg_\calH(h_0)$.  Let $g_0=\Phi_\varphi^{-1}(x_0,x'_0,h_0)$.  Since $\src_\calG$ is locally subductive, after shrinking $U_p$, there is a lift $g_u$ of $x_u$ against $\src_\calG$ through $g_0$.  Then $\varphi(g_u)$ is the desired lift of $y_u$.  Since $\inv_\calH$ is a diffeomorphsm, both $\src_\calH$ and $\trg_\calH$ are local subductions.

		Continuing with the same plot $p\colon u\mapsto y_u$, let $x_0\in\calG_0$ such that $\varphi(x_0)=y_0$. Since $\varphi_0$ is subductive, after shrinking $U_p$, there is a lift $x'_u$ of $p$ to $\calG_0$.  Let $g_0=\Phi_\varphi^{-1}(x_0,x'_0,\unit_{y_0})$.  Since $\trg_\calG$ is locally subductive, after shrinking $U_p$, there is a lift $g_u$ of $x'_u$ against $\trg_\calG$ through $g_0$.  Then $\src_\calG(g_u)$ is the desired lift of $p$, from which it follows that $\varphi_0$ is locally subductive.

		Conversely, suppose $\varphi$ is a subductive weak equivalence with $\src_\calH$, $\trg_H$, and $\varphi_0$ locally subductive.  Let $p\colon u\mapsto x_u$ be a plot of $\calG_0$ and fix $g_0\in\calG_1$ such that $\src_\calG(g_0)=x_0$.  Since $\src_\calH$ is locally subductive, after shrinking $U_p$, there exists a lift $h_u$ of $\varphi(x_u)$ against $\src_\calH$ through $\varphi(g_0)$.  Since $\varphi_0$ is locally subductive, after shrinking $U_p$, there exists a lift $x'_u$ of $\trg_\calH(h_u)$ through $\trg_\calG(g_0)$.  Then $g_u=\Phi_\varphi^{-1}(x_u,x'_u,h_u)$ is the desired lift, from which it follows that $\src_\calG$ and $\trg_\calG$ are locally subductive.  This proves \cref{i:local subd2}.
	\end{proof}
	
	\begin{proposition}[Weak Equivalence between Lie Groupoids]\labell{p:weak equiv Lie}
		Given a weak equivalence $\varphi\colon\calG\to\calH$ between Lie groupoids, $\varphi$ is a Lie weak equivalence. Consequently, $\calW_\Lie=\calW\cap\LieGpoid_1$. Furthermore, if $\varphi$ is also subductive, then it is surjective submersive.  Consequently, $\calJ_\Lie=\calJ\cap\LieGpoid_1$.
	\end{proposition}
	
	\begin{proof}
		Since $\calG_1$ is a manifold and $\Phi_\varphi$ is a diffeomorphism, the codomain of $\Phi_\varphi$ is also a manifold. Since $\trg_\calH$ is a surjective submersion, $\calG_0\ftimes{\varphi}{\trg}\calH_1$ is a manifold.  Since $\calG$ is Lie, $\src_\calG$ and $\trg_\calG$ are surjective submersions and hence locally subductive by \cref{i:submersion} of \cref{x:ind subd}.  By \cref{i:local subd1} of \cref{l:weak equiv local subd}, $\Psi_\varphi$ is a local subduction, and hence a surjective submersion.  Thus $\varphi$ is a Lie weak equivalence.
		
		If $\varphi$ is a subductive weak equivalence, then $\varphi_0$ is locally subductive by \cref{i:local subd2} of \cref{l:weak equiv local subd}, and hence a surjective submersion by \cref{i:submersion} of \cref{x:ind subd}.
	\end{proof}	
	
	A diffeological groupoid being Lie (or \emph{not} being Lie) is not a Morita invariant.  For instance, the trivial groupoid of a point, which is Lie, is Morita equivalent to the pair groupoid of any diffeological space.  However, some diffeological groupoids that are not Lie do \emph{not} admit a Morita equivalence to any Lie groupoid at all:
	
	\begin{example}[$\ZZ/2\circlearrowright\RR$]\labell{x:not local subd}
		Let $\ZZ/2$ act on $\RR$ by reflection.  The corresponding action groupoid $\calG$ admits a kernel $\calK:=\ker(\chi_\calG)$ to its characteristic functor $\chi_\calG:=(\src_\calG,\trg_\calG)\colon\ZZ/2\ltimes\RR\to(\RR^2\toto\RR)$ (see \cref{i:relation gpd,i:kernel gpd} of \cref{x:diffeol gpd}).  The source fibres of $\calK$ are made up solely of the units of $\calG$ except for the source fibre of $0$, which is isomorphic to $\ZZ/2$.  It follows that $\src_\calK$ is not locally subductive, as one cannot lift a non-trivial path through $0$ of the object space $\RR$ to a path through the non-trivial arrow in the stabiliser at $0$.
		
		Given a Morita equivalence $\calK\ot{\varphi}\calL\underset{\raisebox{1mm}{$\scriptstyle\psi$}}{\we}\calH$, the maps $\src_\calL$ and $\trg_\calL$ cannot be locally subductive by \cref{i:local subd2} of \cref{l:weak equiv local subd}.  By \cref{r:morita}, we can always choose $\calL$ and $\psi$ such that $\psi$ is a subductive weak equivalence, in which case \cref{i:local subd2} of \cref{l:weak equiv local subd} also implies that $\src_\calH$ and $\trg_\calH$ cannot be local subductions.  Thus $\calH$ cannot be a Lie groupoid.
	\end{example}
	
	Let $\LieGpoid_\ana$ be the bicategory of Lie groupoids with anafunctors (in which the weak equivalences on the left of each anafunctor is surjective submersive) as $1$-cells and transformations as defined in \cref{d:transformation} as $2$-cells.  This is the localisation of the $2$-site $(\LieGpoid,\calJ_\Lie)$ at $\calW_\Lie$ \cite{roberts:internal,roberts:anafunctors}.  To show that this is an essentially full sub-bicategory of $\DGpoid_\ana$, given \cref{p:weak equiv Lie}, it suffices to show that every anafunctor between two Lie groupoids admits a $2$-cell to a Lie anafunctor. 
	
	\begin{proposition}[Anafunctors Versus Lie Anafunctors]\labell{p:ana-lie}
		Let $\calG$ be a Lie groupoid and $\calG\ot{\varphi}\calK\underset{\psi}{\longrightarrow}\calH$ an anafunctor in $\DGpoid_\ana$.  There is a Lie groupoid $\calK'$, an anafunctor $\calG\ot{\varphi'}\calK'\underset{\psi'}{\longrightarrow}\calH$, and a transformation between the two anafunctors.  In particular, any anafunctor between Lie groupoids admits a transformation to a Lie anafunctor.  Moreover, if $\calG$ is Hausdorff (resp.\ second-countable), then so is $\calK'$.
	\end{proposition}

	\begin{proof}
		Let $\calU=\{i_\mu\colon U_\mu\hookrightarrow\calG_0\}_{\mu\in A}$ be a countable open cover of $\calG_0$.  Without loss of generality, assume that each $U_\mu$ is an open subset of $\RR^n$ where $n=\dim\calG_0$.  Then $\calU$ is a generating family of the diffeology on $\calG_0$.  Since $\varphi_0$ is a subduction, for each $\mu$ there is a countable open cover $\calV_\mu=\{V_{\mu,\nu}\hookrightarrow U_\mu\}_{\nu\in B_\mu}$ of $U_\mu$ such that $i_\mu|_{V_{\mu,\nu}}$ admits a lift $j_{\mu,\nu}$ to $\calK_0$.  Let $\calV:=\{\varphi\circ j_{\mu,\nu}\}_{\mu\in A,\nu\in B_\mu}$.  Then $\calV$ is a generating family of the diffeology on $\calG_0$, $\Neb(\calV)$ is a manifold, and $\ev_\calV\colon\Neb(\calV)\to\calG_0$ is a surjective submersion.  By \cref{x:pullback surj subm}, the pullback $\ev_\calV^*\calG$ is a Lie groupoid, and by \cref{i:subdwe pullback} of \cref{l:we properties}, $\varphi':=(\ev_\calV,\pr_3)$ sending objects $v\in V_{\mu,\nu}$ to $\varphi(j_{\mu,\nu}(v))$ and arrows $(v_1,v_2,g)\in\Neb(\calV)^2\ftimes{\ev_\calV^2}{(\src,\trg)}\calG_1$ to $g$ is a subductive weak equivalence.

		Due to its construction, there is a local diffeomorphism $I$ from $\Neb(\calV)$ to $\Neb(\calD_{\calK_0})$, the nebula of the diffeology $\calD_{\calK_0}$ with evaluation map $\ev_\calK$, defined as follows.  For each $q\in\calV$, there is a plot $p_q$ of $\calK_0$ so that $q=\varphi\circ p_q$, in which case $U_q=U_{p_q}\subseteq\Neb(\calD_{\calK_0})$.  Define $I|_{U_q}:=\id_{U_{p_q}}$ (as a side remark, the Axiom of Choice is required here).  Let $\psi'\colon\ev_\calV^*\calG\to\calH$ be the smooth functor given by $\psi'_0:=\psi\circ\ev_\calK\circ I$ and $$\psi'_1(v_1,v_2,g):=\psi\left(\Phi_\varphi^{-1}(\ev_\calK(v_1),\ev_\calK(v_2),g)\right).$$  Set $\calK'=\ev_\calV^*\calG$ and let $\calL:=\calK\ftimes{\varphi}{\varphi'}\calK'$ and $S\colon\psi\circ\pr_1\Rightarrow\psi'\circ\pr_2$ be given by the smooth map $$S\colon\calL_0\to\calH_1\colon(y,v)\mapsto\psi\left(\Phi^{-1}_{\varphi}(y,\ev_\calK(v),\unit_{\varphi(y)})\right).$$  Then $\calG\ot{\varphi'}\calK'\uto{\psi'}\calH$ is an anafunctor, and it remains to show that $S$ is a natural transformation.

		Fix an arrow $(k,(v_1,v_2,g))$ from $(y_1,v_1)$ to $(y_2,v_2)$ in $\calL_1$.  We need to show that
		$$\psi(\Phi_\varphi^{-1}(\ev_\calK(v_1),\ev_\calK(v_2),g)\cdot\Phi_\varphi^{-1}(y_1,\ev_\calK(v_1),\unit_{\varphi(y_1)})) = \psi(\Phi_\varphi^{-1}(y_2,\ev_\calK(v_2),\unit_{\varphi(y_2)})\cdot k).$$  Since $\varphi(k)=g$ and $\psi$ is a functor, both sides of the equality reduce to 
		$\Phi_\varphi^{-1}(y_1,\ev_\calK(v_2),g)$, proving naturality.
		
		If $\calH$ is also a Lie groupoid, then \cref{p:weak equiv Lie} guarantees that $\calG\ot{\varphi'}\calK'\uto{\psi'}\calH$ is a Lie anafunctor.  Finally, if $\calG$ is Hausdorff (resp.\ second-countable), then so is $\calK'$ by \cref{p:pullback}, completing the proof.
	\end{proof}
	
	Our goal for this section now follows:
	
	\begin{theorem}[$\LieGpoid_\ana$ in $\DGpoid_\ana$]\labell{t:ess full subbicateg}
		The bicategory $\LieGpoid_\ana$ is an essentially full sub-bicategory of $\DGpoid_\ana$.  Moreover, the restriction to Hausdorff and/or second-countable Lie groupoids forms an essentially full sub-bicategory of $\LieGpoid_\ana$, and hence also of $\DGpoid_\ana$.
	\end{theorem}
	
	One may wonder if the inclusion of $\LieGpoid_\ana$ into $\DGpoid_\ana$ is indeed a pseudofunctor, in particular with respect to the various associators, unitors, and compositions.  But these are all defined in precisely the same way in $\LieGpoid_\ana$ as they are in $\DGpoid_\ana$, and so the inclusion preserves these constructions.
	
	Restricting our attention to Morita equivalence:
	
	\begin{corollary}\labell{c:morita equivalence}
		Given a (diffeological) Morita equivalence between two Lie groupoids $\calG\ot{\varphi}\calK\overset{\simeq}{\uto{\psi}}\calH$, there is a Lie anafunctor that is a Morita equivalence $\calG\ot{\varphi'}\calK'\overset{\simeq}{\uto{\psi'}}\calH$ and a transformation between the two anafunctors.  Thus, the Morita equivalence between two Lie groupoids in the diffeological sense is equivalent to that in the Lie sense.
	\end{corollary}

%%%%%%%%%%%%%%%%%%%%%%%%%%%%%%%%%%%%%%%%%%%%%%%%
\section{Diffeological Bibundles}\labell{s:bibundles}
%%%%%%%%%%%%%%%%%%%%%%%%%%%%%%%%%%%%%%%%%%%%%%%%	

	In \cite{vdS:thesis,vdS:morita}, van der Schaaf develops the bicategory of diffeological groupoids with bibundles, along with the notion of Morita equivalence in this context.  He leaves open the following question: Is a diffeological Morita equivalence between Lie groupoids necessarily a Morita equivalence in the Lie sense?  We almost answered the question affirmatively above via \cref{c:morita equivalence}.  The only thing to check is that Morita equivalence in terms of bibundles is the same thing as Morita equivalence in the anafunctor context.  This will be accomplished by showing that the bicategory using bibundles is equivalent to that using anafunctors.
	
	In the following we refer the reader to \cite{vdS:morita} for full details and definitions.

	\begin{definition}[Bibundles]\labell{d:bibundle}
		Given diffeological groupoids $\calG$ and $\calH$, a \textbf{(diffeological) $(\calG,\calH)$-bibundle} comprises a left action $\calG\LAalong{l_X}X$ and a right action $X\RAalong{r_X}\calH$ on a diffeological space $X$ such that the left anchor map $l_X$ is $\calH$-invariant, the right anchor map $r_X$ is $\calG$-invariant, and the actions commute.  Denote these by $\calG\LAalong{l_X}X\RAalong{r_X}\calH$. The bibundle is \textbf{left principal} if the underlying left bundle $\calG\LAalong{l_X}X\overset{r_X}\longrightarrow\calH_0$ is principal; it is \textbf{right principal} if the underlying right bundle $\calG_0\overset{l_X}{\longleftarrow}X\RAalong{r_X}\calH$ is principal. It is \textbf{biprincipal} if it is both left and right principal.  We often represent a bibundle diagrammatically as follows:
		$$\xymatrix{
			\calG_1 \ar@<0.5ex>[d] \ar@<-0.5ex>[d] \ar@{}[r] |{\curvearrowright} & X \ar[dl]^{l_X} \ar[dr]_{r_X} & \calH_1 \ar@<0.5ex>[d] \ar@<-0.5ex>[d] \ar@{}[l] |{\curvearrowleft}\\
			\calG_0 & & \calH_0. \\
		}$$

		For a fixed diffeological groupoid $\calG$, the \textbf{identity bibundle} is given by $\calG\LAalong{\trg_\calG}\calG_1\RAalong{\src_\calG}\calG$, where the actions are as given by left and right groupoid multiplication.  
	\end{definition}
	
	The $1$-cells of the bibundle bicategory will be right-principal bibundles, and so we focus on those.  The composition of bibundles is defined using a construction known as the ``balanced tensor product''; also known as the the Hilsum-Skandalis tensor product in the literature, see \cite{HS}.
	
	\begin{definition}[Balanced Tensor Products]\labell{d:balanced tensor product}
		Given right principal bibundles $\calG\LAalong{l_X}X\RAalong{r_X}\calH$ and $\calH\LAalong{l_Y}Y\RAalong{r_Y}\calK$ define the \textbf{balanced tensor product} to be the space $X\otimes_{\calH}Y:=(X\ftimes{r_X}{l_Y}Y)/\calH$, where the $\calH$-action is the antidiagonal action $((x,y),h)\mapsto(x\cdot h,h^{-1}\cdot y)$.  Equip $X\otimes_\calH Y$ with a left $\calG$-action $$\calG\LAalong{L_X}X\otimes_\calH Y\colon g\cdot(x\otimes y):=(gx)\otimes y$$ with left anchor map $$L_X\colon X\otimes_\calH Y\to \calG_0\colon x\otimes y\mapsto l_X(x),$$ and with a right $\calK$-action $$X\otimes_\calH Y\RAalong{R_Y}\calK\colon (x\otimes y)\cdot k:= x\otimes(yk)$$ with right anchor map $$R_Y\colon X\otimes_\calH Y\to \calK_0\colon x\otimes y \mapsto r_Y(y).$$
This is a well-defined right principal bibundle \cite[Constructions 4.6, 5.8]{vdS:morita}.
	\end{definition}

	\begin{definition}[Bi-equivariant Diffeomorphisms]\labell{d:bieqvt diffeom}
		Given diffeological groupoids $\calG$ and $\calH$, a \textbf{bi-equivariant diffeomorphism} from a right-principal bibundle $\calG\LAalong{l_X}X\RAalong{r_X}\calH$ to a right-principal bibundle $\calG$-bundle $\calG\LAalong{l_Y}Y\RAalong{r_Y}\calH$ is a diffeomorphism $\alpha\colon X\to Y$ that is both $\calG$-equivariant and $\calH$-equivariant (\emph{i.e.}\ \textbf{bi-equivariant}), such that $l_X=l_Y\circ\alpha$ and $r_X=r_Y\circ\alpha$.  Thus, the following diagram commutes:
		$$\begin{gathered}[b]\xymatrix{
			\calG_1 \ar@<0.5ex>[d] \ar@<-0.5ex>[d] \ar@{}[r] |{\curvearrowright} & X \ar[dl]^{l_X} \ar[dr]_{r_X} \ar[dd]_{\alpha} & \calH_1 \ar@<0.5ex>[d] \ar@<-0.5ex>[d] \ar@{}[l] |{\curvearrowleft} \\
			\calG_0 & & \calH_0 \\
			\calG_1 \ar@<0.5ex>[u] \ar@<-0.5ex>[u] \ar@{}[r] |{\rotatebox{180}{$\scriptstyle\curvearrowleft$}} & Y \ar[ul]_{l_Y} \ar[ur]^{r_Y} & \calH_1 \ar@<0.5ex>[u] \ar@<-0.5ex>[u] \ar@{}[l] |{\rotatebox{180}{$\scriptstyle\curvearrowright$}} \\
		}\\[-\dp\strutbox]\end{gathered}\eqno\qedhere$$
	\end{definition}
	
	The $2$-cells of the bibundle bicategory will be the bi-equivariant diffeomorphisms between the bibundles.  
	
	\begin{theorem}[The Bicategory $\DBiBund$]\labell{t:DBiBund}
		There is a bicategory $\DBiBund$ consisting of diffeological groupoids as objects, right principal bibundles as $1$-cells, and bi-equivariant diffeomorphisms as $2$-cells.  
	\end{theorem}	
	
	See \cite{vdS:morita} for details on the unitors, associators, etc. The proof follows from \cite[Theorem 5.17, Subsection 5.3]{vdS:morita}; the theorem there focuses on $1$-cells as bibundles that are not necessarily principal, and $2$-cells that are not necessarily diffeomorphisms (this yields a category for diffeological groupoids, although it does not for Lie groupoids).  Subsection 5.3 then shows that the restriction to the sub-bicategory as in \cref{t:DBiBund} satisfies the required coherence relations and identities.  The proof of \cite[Theorem 5.17]{vdS:morita} is analogous to that for Lie groupoid theory; see \cite[Proposition 2.12]{blohmann:stacky}.

	Similar to spanisation, there is a natural way to turn a smooth functor between diffeological groupoids into a right principal bibundle, and a smooth natural transformation into a bi-equivariant diffeomorphism.

	\begin{definition}[Bibundlisation]\labell{d:bibundlisation}
		Let $\varphi\colon\calG\to\calH$ be a smooth functor.  Its \textbf{bibundlisation} (or just bundlisation, as in \cite{blohmann:stacky}) is the right principal bibundle $$\calG\LAalong{\pr_1}(\calG_0)\ftimes{\varphi}{\trg}\calH_1\RAalong{\Psi_\varphi}\calH,$$ where the $\calG$-action sends $(g,(x,h))$ to $(\trg_\calG(g),\varphi(g)h)$ with anchor map $\pr_1$, and the $\calH$-action sends $((x,h),h')$ to $(x,hh')$.  If $\psi\colon\calG\to\calH$ is another smooth functor and $S\colon\varphi\Rightarrow\psi$ is a smooth natural transformation, then the \textbf{bibundlisation} of $S$ is the bi-equivariant diffeomorphism $\alpha\colon(\calG_0)\ftimes{\varphi}{\trg}\calH_1\to(\calG_0)\ftimes{\psi}{\trg}\calH_1$ sending $(x,h)$ to $(x,S(x)h)$.
	\end{definition}

	\begin{remark}[Bibundlisation versus Anafunctisation of Functors]\labell{r:bibundlisation}
		Given a smooth functor $\varphi\colon\calG\to\calH$, the ($\calG$-$\calH$)-action groupoid corresponding to the joint $\calG$- and $\calH$-actions on $(\calG_0)\ftimes{\varphi}{\trg}\calH_1$, denoted $\calG\ltimes(\calG_0)\ftimes{\varphi}{\trg}\calH_1\rtimes\calH$, is defined to be the action groupoid of the left groupoid action of $\calG\times\calH$ on $(\calG_0)\ftimes{\varphi}{\trg}\calH_1$ with action map $((g',h'),(x,h))\mapsto(\trg_\calG(g'),\varphi(g')h(h')^{-1})$ and anchor map $(x,h)\mapsto(x,s_\calH(h))$.  It is straightforward to check that the map $\calG\ltimes(\calG_0)\ftimes{\varphi}{\trg}\calH_1\rtimes\calH\to\calG\wtimes{\varphi}{\id_\calH}\calH$ sending $((g',h'),(x,h))$ to $(g',h^{-1},h')$ is an isomorphism of diffeological groupoids.  In particular, the action groupoid of the bibundlisation of $\varphi$ is isomorphic to the anafunctisation of $\varphi$; see \cref{r:quasi-inverse}.
	\end{remark}

	\begin{lemma}[Bibundlisation of a Weak Equivalence]\labell{l:bibundlisation}
		The bibundlisation of a smooth functor is biprincipal if and only if the functor is a weak equivalence.
	\end{lemma}

	\begin{proof}
		We follow the terminology of \cite[Section 5]{vdS:morita}.  Let $\varphi\colon\calG\to\calH$ be a smooth functor.  It is immediate from its definition that the right anchor map of the bibundlisation of $\varphi$ is  subductive if and only if $\varphi$ is smoothly essentially surjective.  Moreover, if $\varphi$ is smoothly fully faithful, then the action map $$A_\calG\colon\calG\ftimes{\src}{\pr_1}\left((\calG_0)\ftimes{\varphi}{\trg}\calH_1\right)\to\left((\calG_0)\ftimes{\varphi}{\trg}\calH_1\right)\ftimes{\Psi_\varphi}{\Psi_\varphi}\left((\calG_0)\ftimes{\varphi}{\trg}\calH_1\right)$$ is a diffeomorphism.  Conversely, if $A_\calG$ is a diffeomorphism, then the division map of the $\calG$-bundle $(\calG_0)\ftimes{\varphi}{\trg}\calH_1\underset{\Psi_\varphi}{\longrightarrow}\calH_0$ is well-defined and smooth, from which it follows that $\Phi_\varphi$ is a diffeomorphism. 
	\end{proof}

	\cref{l:bibundlisation} indicates that biprincipal bibundles may be the bibundle version of the anafunctors admitting quasi-inverses as in \cref{p:quasi-inverse}.  This is confirmed in \cite[Proposition 5.24]{vdS:morita}.

	\begin{theorem}\labell{t:bibundle-localise}
		Bibundlisation $\mathfrak{B}\colon\DGpoid\to\DBiBund$, sending objects to themselves and everything else to their bibundlisation, is a localisation of $\calW$.  In particular, $\DGpoid_\ana$ and $\DBiBund$ are equivalent bicategories.
	\end{theorem}

	\begin{proof}
		We begin by showing that bibundlisation $\mathfrak{B}$ is a pseudofunctor.  This is straightforward bookkeeping, and we refer the reader to other sources such as \cite{JY} for all of the bicategorical definitions for the sake of brevity.  By definition, $\mathfrak{B}$ sends diffeological groupoids to themselves, and smooth functors and smooth natural transformations to their bibundlisations (see \cref{d:bibundlisation}).

		Fixing diffeological groupoids $\calG$ and $\calH$, $\mathfrak{B}$ is required to restrict to a functor from $\DGpoid(\calG,\calH)$ to $\DBiBund(\calG,\calH)$.  But this is immediate, since the vertical composition of bi-equivariant diffeomorphisms is just their composition in the standard sense \cite[Proposition 5.14]{vdS:morita}, and the bibundlisation of an identity natural transformation is equal to the identity bi-equivariant diffeomorphism.

		Given a diffeological groupoid $\calG$, there is a bi-equivariant diffeomorphism $\iota_\calG$ from the identity bibundle of $\calG$ to the bibundlisation of $\id_\calG$, given by the inverse of the map $\pr_2\colon(\calG_0)\ftimes{\id_\calG}{\trg}\calG_1\to\calG_1$.  Also, given functors $\varphi\colon\calG\to\calH$ and $\psi\colon\calH\to\calK$, there is a bi-equivariant diffeomorphism $$\gamma_{\varphi,\psi}\colon\left((\calG_0)\ftimes{\varphi}{\trg}\calH_1\right)\otimes_\calH\left((\calH_0)\ftimes{\psi}{\trg}\calK_1\right)\to(\calG_0)\ftimes{\psi\circ\varphi}{\trg}\calK_1\colon(x,h)\otimes(y,k)\mapsto(x,\psi(h)k).$$  We need to confirm the coherence relations for a pseudofunctor; see \cite[(4.1.3) and (4.1.4)]{JY} or \cite[(M.1) and (M.2)]{benabou}.
		
		For the pseudofunctorial associativity coherence relation \cite[(4.1.3)]{JY}, fix three functors $$\calG\underset{\varphi}{\longrightarrow}\calH\underset{\psi}{\longrightarrow}\calK\underset{\chi}{\longrightarrow}\calL.$$  Let $A\colon(x\otimes y)\otimes z\mapsto x\otimes(y\otimes z)$ be the associator of $\DBiBund$ \cite[Proposition 5.13]{vdS:morita}.  Then the coherence relation reduces to showing that $$\gamma_{\psi\circ\varphi,\chi}\circ(\gamma_{\varphi,\psi}\otimes\id_{\calK_0\times\calL_1})=\gamma_{\varphi,\chi\circ\psi}\circ(\id_{\calG_0\times\calH_1}\otimes\gamma_{\psi,\chi})\circ A.$$  The left-hand side sends $((x,h)\otimes(\src_\calH(h),k))\otimes(\src_\calK(k),\ell)$ to $(x,\chi(\psi(h)k)\ell)$, whereas the right-hand side sends the same point to $(x, \chi\circ\psi(h)\chi(k)\ell)$; these are equal.

		Focusing on $\varphi\colon\calG\to\calH$, let $\lambda^{\bi}_{\calG,\calH}$ and $\rho^{\bi}_{\calG,\calH}$ be the left and right unitors of $\DBiBund$ \cite[Proposition 5.12]{vdS:morita}.  For the pseudofunctorial left and right unity coherence relations \cite[(4.1.4)]{JY}, the first reduces to $$\lambda^{\bi}_{\calG,\calH}(\calG_0\ftimes{\varphi}{\trg}\calH_1)=\mathfrak{B}(\id_\varphi)\circ\gamma_{\id_\calG,\varphi}\circ(\iota_\calG\otimes\id_{\calG_0\times\calH_1}).$$  Both sides send $(g,(\src_\calG(g),h))$ to $(\trg_\calG(g),\varphi(g)h)$, confirming this coherence relation.  The second relation reduces to $$\rho^{\bi}_{\calG,\calH}(\calG_0\ftimes{\varphi}{\trg}\calH_1)=\mathfrak{B}(\id_\varphi)\circ\gamma_{\varphi,\id_\calH}\circ(\id_{\calG_0\times\calH_1}\otimes\iota_\calH),$$ in which both sides send $((x,h),h')$ to $(x,hh')$, confirming the last coherence relation.  We conclude that $\mathfrak{B}$ is a pseudofunctor.

		Next, we use Pronk's Comparison Theorem \cite[Proposition 24]{pronk} to show that $\mathfrak{B}$ `extends' to an equivalence of bicategories $\DGpoid_\ana$.  We need to show
	\begin{enumerate}
		\item\labell{i:B sends W to adj equiv} $\mathfrak{B}$ sends weak equivalences to equivalences in $\DBiBund$;
		\item\labell{i:pronk-ess surj bi} $\mathfrak{B}$ is essentially surjective on objects (this is trivially true for $\mathfrak{B}$);
		\item\labell{i:pronk-compatibility bi} for every bibundle $\calH\LAalong{l_Y}Y\RAalong{r_Y}\calK$, there exist a diffeological groupoid $\calG$, a weak equivalence $\varphi\colon\calG\to\calH$, a smooth functor $\omega\colon\calG\to\calK$, and a bi-equivariant diffeomorphism $\alpha$ from $\mathfrak{B}(\omega)$ to the composition of $\calH\LAalong{l_Y}Y\RAalong{r_Y}\calK$ with $\mathfrak{B}(\varphi)$; and
		\item\labell{i:pronk-ff bi} $\mathfrak{B}$ is fully faithful on $2$-cells.
	\end{enumerate}

		Let $\varphi\colon\calG\we\calH$ be a weak equivalence.  By \cref{l:bibundlisation} and \cite[Proposition 5.24]{vdS:morita}, $\mathfrak{B}(\varphi)$ is a biprincipal bibundle, an equivalence.  Thus \cref{i:B sends W to adj equiv} is satisfied.

		For \cref{i:pronk-compatibility bi}, let $\calG$ be the action groupoid $\calH\ltimes Y\rtimes\calK$, define $\varphi\colon((h,k),y)\mapsto h$, and $\omega\colon((h,k),y)=k^{-1}$.  Define $\beta((y,h)\otimes y'):=(y,d_\calK(y',h^{-1}y))$, where $d_\calK$ is the division map of the underlying $\calK$-bundle.  This is well-defined and smooth with smooth inverse $\alpha(y,k)=(y,\unit_{r_Y(y)})\otimes yk$.  Since $\varphi_0=l_Y$ and $\omega_0=r_Y$, it follows that $\beta$ is bi-equivariant.  This proves \cref{i:pronk-compatibility bi}.

		Finally, for \cref{i:pronk-ff bi}, let $\varphi,\psi\colon\calG\to\calH$ be smooth functors with $S_1,S_2\colon\varphi\Rightarrow\psi$ smooth natural transformations.  If $\mathfrak{B}(S_1)=\mathfrak{B}(S_2)$, then by definition, $(x,S_1(x)h)=(x,S_2(x)h)$ for every $(x,h)\in(\calG_0)\ftimes{\varphi}{\trg}\calH_1$, from which $S_1=S_2$ follows.  On the other hand, for any bi-equivariant diffeomorphism $\alpha\colon(\calG_0)\ftimes{\varphi}{\trg}\calH_1\to\calG\ftimes{\psi}{\trg}\calH_1$, define $S(x):=h'h^{-1}$ for any $h,h'\in\calH_1$ such that $(x,h)\in(\calG_0)\ftimes{\varphi}{\trg}\calH_1$ and $(x,h')=\alpha(x,h)\in\calG_0\ftimes{\psi}{\trg}\calH_1$; this is well-defined by the bi-equivariance of $\alpha$.  Moreover, that $S$ is a natural transformation follows from the $\calG$-equivariance of $\alpha$.  We have shown that $\mathfrak{B}$ acts fully faithfully on $2$-cells.  This completes the proof.
	\end{proof}
	
	Combining \cref{t:ess full subbicateg,t:bibundle-localise}, we answer the open problem of van der Schaaf \cite[Question 7.6]{vdS:morita} affirmatively:

	\begin{corollary}\labell{c:bibundle-localise}
		A biprincipal bibundle between two Lie Groupoids in $\DBiBund$ is a biprincipal bibundle in the Lie sense.  That is, a diffeological Morita equivalence between two Lie groupoids via diffeological bibundles is a Lie Morita equivalence via a bibundle in the Lie groupoid sense.
	\end{corollary}
	
	Note that this result is stronger than one may initially realise: since $2$-cells in $\DBiBund$ are diffeomorphisms, any bibundle representing a Morita equivalence between two Lie groupoids in $\DBiBund$ \emph{must} be a manifold.

%%%%%%%%%%%%%%%%%%%%%%%%%%%%%%%%%%%%%%%%%%%%%%%%
\section{Further Applications \& Examples}\labell{s:apps and examples}
%%%%%%%%%%%%%%%%%%%%%%%%%%%%%%%%%%%%%%%%%%%%%%%%

In this section, we consider certain constructions that remain invariant under Morita equivalence.  We start with the orbit space of a diffeological groupoid, as well as the relation groupoid of the induced equivalence relation on the object space.  Next, we consider the inertia groupoid, and show that this is a Morita invariant.  We then consider principal bundles, connecting the category of anafunctors between the trivial groupoid of a diffeological space and an abelian diffeological group to the corresponding diffeological \v{C}ech cohomology.  Finally, we show that Hausdorff Morita equivalent manifolds with singular foliations yield (diffeologically) Morita equivalent holonomy groupoids.

%% %% %%
\subsection{Orbit Spaces \& Relation Groupoids}\labell{ss:relation}
%% %% %%

	In this subsection, we show that the orbit space and relation groupoid induced by a diffeological groupoid are Morita invariants.  In fact, the relation groupoid comes with two different diffeologies: the first is induced by the pair groupoid as in \cref{i:relation gpd}.  The second is the ``pushforward diffeology'' on the same underlying groupoid induced by the characteristic functor.  In the Lie case, the difference between the two diffeologies gives the obstruction to the groupoid representing a gerbe; see \cite[Definition 4.18]{WW}.
	
	\begin{definition}[Orbit Space]\labell{d:orbit space}
		Given a diffeological groupoid $\calG$, the \textbf{orbit space of $\calG$}, denoted $\calG_0/\calG_1$, is the quotient diffeological space induced by the equivalence relation on $\calG_0$ in which $x\sim x'$ if there exists $g\in\calG_1$ such that $\src_\calG(g)=x$ and $\trg_\calG(g)=x'$.  Denote the quotient map by $\pi_\calG$.
	\end{definition}
	
	We show that the orbit space is a Morita invariant of diffeological groupoids, extending the result for Lie groupoids (\cite[Theorem 3.8]{watts:lgpd-derham}; see also \cite[Theorem 4.3.1]{delhoyo:orbispaces} for a similar result).  In fact, it will be evident from the proof that there is a pseudofunctor from $\DGpoid$ to $\Diffeol$, treating the latter as a trivial $2$-category.
	
	\begin{proposition}[The Orbit Space is a Morita Invariant]\labell{p:orbit space}
		Morita equivalent diffeological groupoids have diffeomorphic orbit spaces.
	\end{proposition}
	
	\begin{proof}
		Let $\calG$ and $\calH$ be Morita equivalent diffeological groupoids.  By definition of Morita equivalence, it suffices to assume that there is a weak equivalence $\varphi\colon\calG\we\calH$.  Define $\widecheck{\varphi}\colon\calG_0/\calG_1\to\calH_0/\calH_1$ by $\widecheck{\varphi}([x]):=[\varphi(x)]$, where the square brackets indicate the image under the quotient map to the orbit space.  This is well-defined by the functoriality of $\varphi$.  Let $p$ be a plot of $\calG_0/\calG_1$.  By definition of the quotient diffeology, after shrinking $U_p$, there is a lift $q$ of $p$ to $\calG_0$.  Then $$\widecheck{\varphi}\circ p=\pi_\calH\circ\varphi\circ q$$ which is smooth.  Smoothness of $\widecheck{\varphi}$ follows.  
		
		Suppose $\widecheck{\varphi}[x]=\widecheck{\varphi}[x']$.  There is an arrow $h\in\calH_1$ such that $\src_\calH(h)=\varphi(x)$ and $\trg_\calH(h)=\varphi(x')$.  Since $\varphi$ is smooth fully faithful, $g:=\Phi_\varphi^{-1}(x,x',h)$ is well-defined, and thus $\src_\calG(g)=x$ and $\trg_\calG(g)=x'$.  Injectivity of $\widecheck{\varphi}$ follows.
		
		Let $p$ be a plot of $\calH_0/\calH_1$.  After shrinking $U_p$, there is a lift $q\colon u\mapsto y_u$ of $p$ to $\calH_0$.  Since $\varphi$ is smoothly essentially surjective, after shrinking $U_p$ again, there is a lift $r\colon u\mapsto (x_u,h_u)$ of $q$ against $\Psi_\varphi$ to $\calG_0\ftimes{\varphi}{\trg}\calH_1$.  Then $$\widecheck{\varphi}(\pi_\calG(x_u))=\pi_\calH(\varphi(x_u))=\pi_\calH(y_u)=p(u).$$  Thus $\widecheck{\varphi}$ is a subduction, and hence a diffeomorphism.
	\end{proof}
	
	The orbit space of a diffeological groupoid is itself Morita equivalent to the relation groupoid induced by the equivalent relation $\sim$ on the object space; see \cref{i:relation gpd} of \cref{x:diffeol gpd}. 
	
	\begin{proposition}[Relation Groupoids]\labell{p:relation gpd}
		Let $X$ be a diffeological space with an equivalence relation $\sim$ and corresponding quotient map $\pi\colon X\to X/\!\sim$  The relation groupoid $X\ftimes{\pi}{\pi}X\toto X$ is weakly equivalent to the trivial groupoid $X/\!\sim$.
	\end{proposition}
	
	\begin{proof}
		Define $\varphi\colon X\ftimes{\pi}{\pi}X\to X/\!\sim$ to be either projection map composed with the quotient map to the quotient.  Then $\varphi_0=\pi$, which is subductive by definition of the quotient diffeology.  
		
		$\Phi_\varphi$ is injective.  Moreover, for any plot $u\mapsto(x_u,x'_u,[x_u])$ of $X^2\ftimes{\pi^2}{\id}(X/\!\sim)$, this lifts to $(x_u,x'_u)$ in $X\ftimes{\pi}{\pi}X$.  Thus $\Phi_\varphi$ is subductive, and hence a diffeomorphism.  Thus $\varphi$ is a weak equivalence.
	\end{proof}
	
	Since nebulaic groupoids of generating families of a diffeology (\cref{i:neb gpd} of \cref{x:diffeol gpd}) are in fact relation groupoids, as the evaluation maps are quotient maps, it follows that these are always Morita equivalent to each other, as well as to the trivial groupoid of the diffeological space itself.
	
	\begin{corollary}[Nebulaic Groupoids]\labell{c:nebulaic gpd}
		Given a diffeological space $X$ and any generating family $\calF$ of the diffeology of $X$, the nebulaic groupoid $\calN(\calF)$ is Morita equivalent to $X$.
	\end{corollary}
	
	We can combine some of the results above, showing that the relation groupoid $(\calG_0)\ftimes{\pi_\calG}{\pi_\calG}(\calG_0)$ induced by a diffeological groupoid $\calG$ is a Morita invariant.
	
	\begin{corollary}[Relation Groupoids of Diffeological Groupoids]\labell{c:relation gpd}
		Given Morita equivalent diffeological groupoids $\calG$ and $\calH$, their corresponding relation groupoids $(\calG_0)\ftimes{\pi_\calG}{\pi_\calG}(\calG_0)$ and $(\calH_0)\ftimes{\pi_\calH}{\pi_\calH}(\calH_0)$ are Morita equivalent.
	\end{corollary}

	We now turn to the image of the characteristic functor $\chi_\calG=(\src_\calG,\trg_\calG)$ of a diffeological groupoid $\calG$.  This has the underlying set-theoretic groupoid of the relation groupoid $(\calG_0)\ftimes{\pi_\calG}{\pi_\calG}(\calG_0)$ as its image; however, the diffeology has fewer plots (it is ``finer'' in the language of diffeology).  This said, it is still a diffeological groupoid.
	
	\begin{proposition}[Images of Characteristic Functors]\labell{p:characteristic functors}
		Given Morita equivalent diffeological groupoids $\calG$ and $\calH$, the images of their corresponding characteristic functors $\chi_\calG(\calG)$ and $\chi_\calH(\calH)$ are Morita equivalent.
	\end{proposition}
	
	\begin{proof}
		It suffices to assume that there is a weak equivalence $\varphi\colon\calG\we\calH$.  We claim that $\varphi_0^2\colon\chi_\calG(\calG)\to\chi_\calH(\calH)$ is a weak equivalence.  It is well-defined and smooth since $\varphi$ is a smooth functor.  
		
		Fix a plot $p\colon u\mapsto y_u$ of $\calH_0$.  Since $\varphi$ is smoothly essentially surjective, after shrinking $U_p$, there is a lift $u\mapsto (x_u,h_u)$ of $p$ to $(\calG_0)\ftimes{\varphi}{\trg}\calH_1$.  The plot $u\mapsto(x_u,(\src_\calH(h_u),\trg_\calH(h_u)))$ is the desired lift of $p$, from which it follows that $\varphi_0^2$ is smoothly essentially surjective.  Injectivity of $\Phi_{\varphi_0^2}$ is immediate from its definition, and subductivity follows from the functoriality of $\varphi$.  The result follows.
	\end{proof}

%% %% %%
\subsection{Inertia Groupoids}\labell{ss:inertia}
%% %% %%
	
We next show that the inertia groupoid of a diffeological groupoid is a Morita invariant; see \cref{i:inertia gpd} of \cref{x:diffeol gpd}.  This is a generalisation of the idea that Morita equivalence should ``preserve the stabilisers'' of diffeological groupoids.  Note that even in the case of Lie groupoids, the inertia groupoids are typically not Lie.

	\begin{proposition}\labell{p:inertia groupoid}
		Given Morita equivalent diffeological groupoids $\calG$ and $\calH$, the inertia groupoids $\calI_\calG$ and $\calI_\calH$ are Morita equivalent.
	\end{proposition}

	\begin{proof}
		It suffices to assume that there is a weak equivalence $\varphi\colon\calG\to\calH$.  Recalling that $\calI_\calG=\calG\ltimes\ker(\chi_\calG)_1$, define $\kappa\colon\calI_\calG\to\calI_\calH$ by $\kappa(g',g):=(\varphi(g'),\varphi(g))$.  This is a well-defined smooth functor since $\varphi$ is.  Let $p\colon u\mapsto h_u$ be a plot of $\ker(\chi_\calH)_1$.  Since $\varphi$ is a weak equivalence, after shrinking $U_p$, there is lift $u\mapsto(x_u,h'_u)$ of $\src_\calH\circ p$ against $\Psi_\varphi$ to $\calG_0\ftimes{\varphi}{\trg}\calH_1$.  Since $h'_uh_u(h'_u)^{-1}$ has source and target $\varphi(x_u)$ for each $u$, $$g_u:=\Phi_\varphi^{-1}(x_u,x_u,h'_uh_u(h'_u)^{-1})$$ is well-defined, smooth, and is a plot of $\ker(\chi_\calG)_1$.  The plot $u\mapsto(g_u,(h'_u,h_u))$ is the desired lift of $p$ against $\Psi_\kappa$, showing that the map is a subduction.
		
		It follows from the injectivity of $\Phi_\varphi$ that $\Phi_\kappa$ is injective.  Fix a plot $p\colon u\mapsto(g_u,\widetilde{g}_u,(h'_u,h_u))$ of $(\calI_\calG)_0^2\ftimes{\kappa^2}{(\src,\trg)}(\calI_\calH)_1$.  Then $h_u=\varphi(g_u)$ and $h'_uh_u(h'_u)^{-1}=\varphi(\widetilde{g}_u)$.  Set $g'_u:=\Phi_\varphi^{-1}(\src_\calG(g_u),\trg_\calG(\widetilde{g}_u),h'_u)$.  Then $\Phi_\kappa(g'_u,g_u)=p(u)$.  Thus $\Phi_\kappa$ is subductive, hence a diffeomorphism.
	\end{proof}

%% %% %% %%
\subsection{Principal Bundles}\labell{ss:princ bdles}
%% %% %% %%

Let $X$ be a diffeological space and $\calG$ a diffeological groupoid.  Consider the category of of all anafunctors from the trivial groupoid of $X$ to $\calG$.  By \cref{t:bibundle-localise}, this category is equivalent to the category of all right principal bibundles between the two diffeological groupoids.  But these bibundles take on a very specific form.

	\begin{lemma}\labell{l:bibundle to princ bdle}
		Let $X$ be a diffeological space and $\calG$ a diffeological groupoid.  The groupoid of bibundles from $X$ to $\calG$ with bi-equivariant diffeomorphisms as arrows is isomorphic to the groupoid of right principal $G$-bundles of $X$ with bundle isomorphisms as arrows.
	\end{lemma}

	\begin{proof}
		Let $X\LAalong{l_Z}Z\RAalong{r_Z}\calG$ be a right principal bibundle.  That $l_Z\colon Z\to X$ is a right principal $\calG$-bundle follows immediately from the definition.  Conversely, if $\rho\colon P\to X$ is a right principal $\calG$-bundle with anchor map $a$, then there is a trivial action of the trivial groupoid of $X$ on $P$, and we have a right principal bibundle $X\LAalong{\rho}P\RAalong{a}\calG$.

		It also follows from the definition of bi-equivariant diffeomorphisms that these are exactly the bundle isomorphisms between the right principal $\calG$-bundles.  The result follows.
	\end{proof}

	If we specialise to the case of a diffeological group $\calG=(G\toto\RR^0)$, then these right principal $\calG$-bundles are exactly the right principal $G$-bundles (see \cite[Definition 5.1]{KWW}).  If we specialise further to an abelian diffeological group $G$, then a result of \cite{KWW} is that principal $G$-bundles over $X$ are classified by the first \v{C}ech cohomology group $\check{H}^1(X;G)$.  In fact, if $\calF$ is a generating family of the diffeology of $X$, then \cite[Lemma 5.11, Remark 5.12, Corollary 5.16]{KWW} imply that the category of right principal bibundles from $\calN(\calF)$ to $\calG$ is equivalent to the category $\check{\calH}^1(\calF,G)$ whose objects are \v{C}ech $1$-cocycles in $\check{C}^1(\calF,G)$, and arrows from $1$-cocycle $f_1$ to $1$-cocycle $f_2$ are the $0$-cochains $\alpha$ in $\check{C}^0(\calF,G)$, such that $\del\alpha=f_2-f_1$.  It now follows from \cref{c:nebulaic gpd} that:

	\begin{theorem}\labell{t:cech}
		The groupoid of anafunctors between $X$ and the abelian group $G$ (viewed as groupoids) is equivalent to $\check{\calH}^1(X,G)$ whose objects are $1$-cocycles in $\check{C}^1(X,G)$ and arrows are $0$-cochains as described above. 
	\end{theorem}

	It now makes sense to explore the groupoid of anafunctors with $2$-cells between two fixed diffeological groupoids $\calG$ and $\calH$ as a generalisation of the \v{C}ech cohomology to $\calG$-equivariant $\calH$-bundles.  But this is outside the scope of this work.

%% %% %% %%
\subsection{Singular Foliations}\labell{ss:fol}
%% %% %% %%

	In \cite{GZ}, the authors define an equivalence relation called Hausdorff Morita equivalence among manifolds equipped with singular foliations.  They then prove that Hausdorff Morita equivalent foliated manifolds yield topologically Morita equivalent holonomy groupoids.  These groupoids are not generally Lie groupoids, as their arrow spaces are quotients of manifolds of varying dimensions.  However, they are naturally diffeological groupoids, as already seen in \cite{AZ}.  We describe the construction of the holonomy groupoid from a diffeological perspective here, and indicate that the main theorem of \cite{GZ} can be promoted to a statement about diffeological Morita equivalence between these holonomy groupoids.

	We refer the reader to \cite{GZ} for definitions of singular foliation, bisubmersion, morphism of bisubmersions, an atlas of bisubmersions of a foliated manifold, and a path holonomy atlas of a foliated manifold.  Let $(M,\calF)$ be a singular foliation of a manifold, and fix an atlas $\calU=(U_i,t_i,s_i)$ of $(M,\calF)$.  Let $G(\calU)$ be the quotient diffeological space $\left(\coprod_{i}U_i\right)/\!\sim$ where for $u\in U_i$ and $v\in U_j$, we have $u\sim v$ if there is an open neighbourhood $U'_i$ of $u$ and a morphism of bisubmersions $f\colon U'_i\to U_j$ sending $u$ to $v$.  Let $Q\colon\coprod_i U_i\to G(\calU)$ be the quotient map.  Define smooth maps $s,t\colon G(\calU)\to M$ by $$s([u]):=s_i(u) \quad \text{and} \quad t([u]):=t_i(u)$$ for $u\in U_i$.  These will serve as source and target maps for a diffeological groupoid structure on $G(\calU)$, with multiplication given by $[u][v]=[f(u,v)]$ where $(u,v)\in U_i\ftimes{s_i}{t_j}U_j$, and $f\colon W\to U_k$ is a morphism of bisubmersions from an open neighbourhood $W$ of $(u,v)$.  Given $x\in M$, the unit of $x$ is given by $[f(u,u)]$ where $u\in U_i$ such that $s_i(u)=x$, and $f\colon W\to U_k$ is a morphism of bisubmersions from an open neighbourhood $W$ of $(u,u)$ in $U_i\ftimes{s_i}{s_i}U'_i$, and $U'_i$ is an open neighbourhood of $u$ in $U_i$ admitting a morphism of bisubmersions $g$ from the inverse bisubmersion $(U'_i,s_i|_{U'_i},t_i|_{U'_i})$ to some element of $\calU$.  Finally, given $[u]\in G(\calU)$, it follows from above that $[g(u)]$ is the inverse of $[u]$.  That all of this is well-defined follows from the definition of an atlas of a foliated manifold, as well as \cite[Corollary 3.4]{GZ}.  It is now straightforward to check that all of these structure maps on $G(\calU)$ are smooth, keeping in mind the fact that each $s_i$ and $t_i$ are submersions.

	Above, if $\calU$ is a path holonomy atlas, then $G(\calU)$ is called the holonomy groupoid of $(M,\calF)$, denoted $H(\calF)$, and it follows from \cite[Corollary 3.17]{GZ} that it is unique up to an isomorphism of diffeological groupoids; here, one can check that the corollary yields diffeologically smooth maps, not just continuous ones.   Finally, \cite[Theorem 3.21]{GZ} can be promoted from an isomorphism between topological groupoids to an isomorphism between diffeological groupoids by straightforward checking. 

	On the other hand, given two singular foliations $(M,\calF_M)$ and $(N,\calF_N)$, these are Hausdorff Morita equivalent if there exists a manifold $P$ with two surjective submersions $\lambda\colon P\to M$ and $\rho\colon P\to N$, each with connected fibres, such that the corresponding pullback foliations $\lambda^{-1}(\calF_M)$ and $\rho^{-1}(\calF_N)$ are equal. 
	
	\begin{theorem}[Morita Equivalent Holonomy Groupoids]\labell{t:gz}
		Given Hausdorff Morita equivalent foliations $(M,\calF_M)$ and $(N,\calF_N)$, the holonomy groupoids $H(\calF_M)$ and $H(\calF_N)$ are Morita equivalent (as diffeological groupoids).
	\end{theorem}

	\begin{proof}
		Let $P$ be a manifold with two surjective submersions with connected fibres $\lambda$ and $\rho$ to $M$ such that $\lambda^{-1}(\calF_M)\cong\rho^{-1}(\calF_N)$.  It follows from the diffeological version of \cite[Theorem 3.21]{GZ} mentioned above that $$\lambda^*H(\calF_M)\cong H(\lambda^{-1}(\calF_M))\cong H(\rho^{-1}(\calF_N))\cong\rho^*H(\calF_N);$$ denote this isomorphism by $\varphi$.  Let $\widehat{\lambda}\colon\lambda^*H(\calF_M)\to H(\calF_M)$ and $\widehat{\rho}\colon\rho^*H(\calF_N)\to H(\calF_N)$ be the corresponding smooth functors.  Then $$H(\calF_M)\ot{\widehat{\lambda}}\lambda^*H(\calF_M)\underset{\widehat{\rho}\circ\varphi}{\subdwe}H(\calF_N)$$ is a Morita equivalent anafunctor.
	\end{proof}
	
	In light of \cref{c:morita equivalence}, if the foliations $(M,\calF_M)$ and $(N,\calF_N)$ in \cref{t:gz} are projective with Hausdorff holonomy groupoids, then the statement becomes an if and only if: this follows from \cite[Proposition 3.39]{GZ} since a diffeological Morita equivalence between Lie groupoids is a Lie Morita equivalence.  It would be interesting to know how much of this theory could be generalised to the theory of integrations of singular subalgebroids of Lie algebroids as in \cite{AZ}.
	
%%%%%%%%%%%%%%%%%%%%%%%%%%%%%%%%%%%%%%%%%%%%%%%%

\printbibliography

\end{document}